\documentclass[
               DIV=15,
               ]{scrartcl}  

\usepackage{color} 
\usepackage{psfrag} 
\usepackage{subfig}
\usepackage{hyperref}
\usepackage{cite}
\usepackage{amsmath}
\usepackage{amssymb}
\usepackage{stmaryrd}

\usepackage{enumitem}
\setlist{nolistsep}
\setlist{noitemsep}

\usepackage{xfrac}
\usepackage{float}
\usepackage{graphics}
\usepackage{grffile}
 \usepackage{pstool} 
\usepackage[textsize=tiny]{todonotes}
\usepackage{placeins}

\usepackage[ruled,linesnumbered]{algorithm2e}
\usepackage{multirow}
\usepackage{stackengine}

\usepackage[square,comma,numbers,sort&compress]{natbib}
\usepackage[export]{adjustbox}

\usepackage{pgfplots}

\usepackage{pgfplotstable}
\usepackage{filecontents}

\usepackage{tikz}
\usepackage{tikzscale}
\usepackage{tikz-3dplot}

\usetikzlibrary{3d}
\usetikzlibrary{spy}
\usetikzlibrary{arrows,shapes}    
\usetikzlibrary{spy}
\usetikzlibrary{backgrounds}  
\usetikzlibrary{decorations}
\usetikzlibrary{decorations.markings}
\usetikzlibrary{positioning}
\usetikzlibrary{patterns}
\usetikzlibrary{calc} 
\usetikzlibrary{fadings}
\usetikzlibrary{trees}
\usetikzlibrary{hobby}
\usetikzlibrary{matrix}
\usepgfplotslibrary{external} 

\usepackage{booktabs}   
\usepackage{tabularx}
\usepackage{diagbox}
\newcolumntype{Y}{>{\centering\arraybackslash}X}

\usepackage{makecell} 
\usepackage{colortbl}   

\usepackage{xspace}
\usepackage{color}     
\usepackage{setspace}   
 
\usepackage{soul}
\usepackage{ulem}
\usepackage{currfile}
\usepackage{import}

\usepackage{csquotes}

\usepackage{authoraftertitle} 
\usepackage{authblk} 

\let\OLDthebibliography\thebibliography
\renewcommand\thebibliography[1]{
  \OLDthebibliography{#1}
  \setlength{\parskip}{0pt}
  \setlength{\itemsep}{0pt plus 0.3ex}
}

\pgfplotsset{compat=1.8}

\newcommand{\findmax}[3]{
    \pgfplotstablesort[sort key={#2},sort cmp={float >}]{\sorted}{#1}%
    \pgfplotstablegetelem{0}{#2}\of{\sorted}%
    \let #3=\pgfplotsretval%
}

\pgfplotscreateplotcyclelist{myCycleListBlackAndWhite}
{%
  {black, mark=o},
  {black, mark=square},
  {black, mark=triangle},
  {black, mark=diamond}, 
  {black, mark=pentagon}, 
  {black, mark=*},
  {black, mark=square*},
  {black, mark=triangle*},
  {black, mark=diamond*}, 
  {black, mark=pentagon*}, 
}

\definecolor{darkgreen}{rgb}{0,0.4,0} 
\definecolor{darkbrown}{rgb}{0.5, 0.396, 0.09}

\pgfplotscreateplotcyclelist{myCycleListColor}
{%
  {blue, mark=o},
  {red, mark=square},
  {darkbrown, mark=triangle},
  {darkgreen, mark=diamond},
  {black, mark=pentagon},
  {blue, mark=*},
  {red, mark=square*},
  {darkbrown, mark=triangle*},
  {darkgreen, mark=diamond*},
  {black, mark=pentagon*},
}       

\pgfplotsset{every axis/.append style= 
              {
                font=\small,
                mark size=2,
                line width = 0.1,
                legend style={font=\small, mark size=3, draw=none, fill=none},
                legend cell align=left,
                cycle list name=myCycleListColor,
              }
            } 
            
\pgfplotstableset
{
  every odd row/.style={before row={\rowcolor[gray]{0.9}}},
  every head row/.style=
  {
    before row=
    {%
      \toprule
    },
    after row=\midrule
  },
  every last row/.style=
  {
    after row=\bottomrule
  }
}

\newif\ifdrawboundingbox

\usetikzlibrary{external} 
\tikzset{external/system call={pdflatex \tikzexternalcheckshellescape
-halt-on-error -interaction=batchmode -jobname "\image" "\texsource"}} 
  
\newcommand{\padnum}[2]{%
  \ifnum#1>1 \ifnum#2<10 0\fi
  \ifnum#1>2 \ifnum#2<100 0\fi
  \ifnum#1>3 \ifnum#2<1000 0\fi
  \ifnum#1>4 \ifnum#2<10000 0\fi
  \ifnum#1>5 \ifnum#2<100000 0\fi
  \ifnum#1>6 \ifnum#2<1000000 0\fi
  \ifnum#1>7 \ifnum#2<10000000 0\fi
  \ifnum#1>8 \ifnum#2<100000000 0\fi
  \ifnum#1>9 \ifnum#2<1000000000 0\fi
  \fi\fi\fi\fi\fi\fi\fi\fi\fi
  \expandafter\expandafter\number#2%
}

\newcounter{tmp}
\newcommand{\nextfigurename}[1]%
{%
	\setcounter{tmp}{\thefigure}%
	\addtocounter{tmp}{1}%
	\tikzsetnextfilename{section\thesubsection_#1}%
}%

\newcommand{\nextsubfigurename}[1]%
{%
	\tikzsetnextfilename{section\thesubsection_#1}%
}%
          
\makeatletter
\renewcommand{\todo}[2][]{\tikzexternaldisable\@todo[#1]{#2}\tikzexternalenable}
\makeatother

\makeatletter
\renewcommand{\missingfigure}[2][]{\tikzexternaldisable\@missingfigure[#1]{#2}\tikzexternalenable}
\makeatother

\newcolumntype{C}[1]{>{\centering\arraybackslash}m{#1}}
\newcolumntype{R}[1]{>{\raggedright\arraybackslash}m{#1}}
\newcolumntype{L}[1]{>{\raggedleft\arraybackslash}m{#1}}

\newcommand{\delete}[1]{\xspace} 
   
\newcommand{\graphDir}{}

\newcommand{\picsDir}{}

\newtheorem{theorem}{Theorem}
\newtheorem{remark}[theorem]{Remark}

\title{Hierarchical multigrid approaches \\ for the finite cell method on uniform \\and multi-level $hp$-refined grids}

\author[1]{J. Jomo\thanks{Corresponding author: john.jomo@tum.de}}
\author[1]{O. Oztoprak}
\author[2]{F. de Prenter}
\author[1]{\authorcr N. Zander}
\author[1]{S. Kollmannsberger}
\author[1]{E. Rank}

\affil[1]{Chair of Computational Modeling and Simulation, Technische Universit\"at M\"unchen}
\affil[2]{Research Develpoment Netherlands, Hengelo, Netherlands}

\usepackage{fancyhdr}
\date{}
\fancypagestyle{plain}{%
  \fancyhf{}%
  \fancyfoot[L]{
  \vspace{-1.5cm} 
  \footnotesize{Preprint submitted to Computer Methods in Applied Mechanics and Engineering}  \hfill
  }}

\drawboundingboxtrue    
\drawboundingboxfalse     

\begin{document}     

\normalem          
\maketitle             
    
\vspace{-1.5cm} 
\hrule 
\section*{Abstract}
This contribution presents a hierarchical multigrid approach for the solution of large-scale finite cell problems on both uniform grids and multi-level $hp$-discretizations. The proposed scheme leverages the hierarchical nature of the basis functions utilized in the finite cell method and the multi-level $hp$-method, which is attributed to the use of high-order integrated Legendre basis functions and overlay meshes, to yield a simple and elegant multigrid scheme. \textcolor{black}{This simplicity is reflected in the fact that all restriction and prolongation operators reduce to binary matrices that do not need to be explicitly constructed}. The coarse spaces are constructed over the different polynomial orders and refinement levels of the immersed discretization. Elementwise and patchwise additive Schwarz smoothing techniques are used to mitigate the influence of the cut cells leading to convergence rates that are independent of the cut configuration, mesh size and in certain scenarios even the polynomial order. \textcolor{black}{The multigrid approach is applied to second-order problems arising from the Poisson equation and linear elasticity. A series of numerical examples demonstrate the applicability of the scheme for solving large immersed systems with multiple millions and even billions of unknowns on massively parallel machines.}

\vspace{0.25cm}
\noindent \textit{Keywords:} immersed methods, finite cell method, multigrid, iterative solvers, \emph{hp}-refinement, parallel computing
 
\vspace{0.25cm}
\hrule  

 
\section{Introduction}
Immersed finite element methods allow numerical simulations to be easily performed on domains with a complex shape by employing a discretization that does not need to capture the boundaries of the body under analysis. The suitability of these methods to handle a variety of geometrical input-data, such as CAD files, scan images, implicitly-defined geometries and oriented point clouds, has led to the emergence of several immersed schemes. Noteworthy immersed methods include the finite cell method (FCM) \cite{Parvizian2007,Duster2008}, cutFEM \cite{Burman2012,Burman2014}, the aggregated unfitted method (AgFEM) \cite{Badia2018a}, the Cartesian grid finite element method \cite{Soriano2013,Navarro2018} and the shifted boundary method \cite{main2018,main2018b}.

The widespread use of immersed methods for solving problems of engineering relevance was for a long time inhibited by the conditioning problems associated with cut elements, i.e.\ elements that are intersected by the boundary of the original body. In \cite{Prenter2017} the underlying cause of ill-conditioning in immersed methods is systematically analyzed for uniform grids and shown to be the occurrence of basis functions that are not only small but also almost linear dependent. Over the past few years, different approaches have been devised to address the conditioning problems of immersed methods such as different preconditioning strategies e.g.\ \cite{Badia2017,Prenter2019,Jomo2019} and methods based on basis function manipulation, removal or aggregation, e.g.\ \cite{Burman2010,Elfverson2018,Badia2018a,Marussig2017review}. These techniques enable iterative solvers to be applied to immersed systems and have opened the door for large-scale immersed finite element analysis.  

One possibility of further improving the convergence properties of (immersed) linear systems and the efficiency of the solution process itself is by the use of multigrid solution techniques. The central idea behind multigrid techniques is to accelerate the convergence of a fine problem by incorporating correction terms generated from a hierarchy
of coarser problems. These terms approximate the smooth components of the error, hereby leading to fast convergence of smooth modes that generally converge slowly when standard iterative solvers such as the Conjugate Gradient (CG) method are used. The high-frequency errors are commonly treated in a smoothing process within the multigrid algorithm. There is a vast amount of literature on multigrid methods. Classical works with a more theoretical perspective include \cite{Hackbusch1982,hackbusch2013}, while \cite{trottenberg2001,Briggs2000,wesseling2004} provide a more implementation-oriented introduction into the subject. Multigrid techniques are generally classified as either geometric or algebraic methods. Geometric multigrid methods generate coarse problems using geometrical information. In finite element analysis, this translates to the generation of coarse problems based on a sequence of discretizations with varying element sizes ($h$-multigrid), polynomial orders ($p$-multigrid) or both ($hp$-multigrid). Algebraic multigrid methods generate their coarse problems solely from the fine system and are usually applied in scenarios where geometrical information is not available or cannot be easily obtained \cite{shapira2003}.

Multigrid methods have been successfully applied to boundary-conforming finite element methods based on the $p$-version of the finite element method e.g.\ \cite{Graig1985,Babuska1989} and isogeometric analysis e.g.\ \cite{Hofreither2016,Gahalaut2013,Tielen2020}, yielding convergence rates that are independent of the grid-size $h$. Convergence rates that are independent of the polynomial order $p$ are, however, harder to achieve and are dependent on the type of smoother employed, as well as the problem type. In most cases, standard smoothers such as the Jacobi or Gauss-Seidel methods result in convergence rates that deteriorate with increasing values of $p$, as shown in \cite{Gahalaut2013}. In \cite{Tielen2020}, it is shown that using an ILUT smoother in an isogeometric $p$-multigrid framework results in $p$-independent convergence rates for the problem classes considered therein. Likewise, smoothers based on the multiplicative Schwarz algorithm have been shown to yield convergence rates independent of $p$ for two-dimensional problems in isogeometric analysis \cite{Riva2019}. The use of different smoothers in the context of high-order isogeometric analysis is discussed in \cite{Donatelli2015,Donatelli2017,Bracco2019}. 

Approaches based on multigrid cycles have also been utilized to solve linear systems derived from immersed finite element methods. An algebraic multigrid solver is utilized in \citep{Verdugo2019} to solve large problems with up to 300 million unknowns using AgFEM on meshes consisting of linear finite elements. It is also used in \cite{badia2020aggregated,badia2020generic} to solve systems with up to 480 million unknowns involving local $h$-refinements using an adaptive $h$-AgFEM framework. Reference \cite{Nussig2018} employs an algebraic multigrid method for the solution of cutFEM discretizations of the electroencephalography (EEG) forward problem. In \cite{Prenter2019b} an $h$-multigrid approach that employs a multiplicative Schwarz smoother is presented and convergence rates independent of the cut configuration and the element size are reported for Lagrange and B-spline bases in problems with up to 10 million unknowns. Reference \cite{Saberi2020} employs a geometric multigrid approach to solve finite cell systems arising from the Poisson equation on a square domain. This paper shows promising developments for the solution of large FCM systems but only reports results of two-dimensional problems with linear elements on which local $h$-refinement is performed.  
 
As previously mentioned, it is possible to construct the nested spaces of the multigrid algorithm in high-order and $hp$-adaptive finite element methods based on the values of the polynomial order $p$ and the different refinement levels of the mesh. In a $p$-multigrid, the nested subspace of a certain polynomial order is simply spanned by the basis functions up to and including this order. $p$-multigrid methods were first developed in the context of the $p$-version of the finite element method \cite{Graig1985} and the spectral element method \cite{Roquist1987}. Since then $p$-multigrid approaches have been extensively studied, see e.g.\ \citep{Babuska1989,Yserentant1985,Graig1985,Yserentant1986,Foresti1989} and applied to different problem classes. Most of these methods define the multigrid spaces by using an arithmetic sequence where the polynomial order is successively lowered until $p=1$, see \cite{Mitchell2010}. It is, however, also possible to define the V-cycles by using a geometrical sequence of $p$, as suggested in \cite{Babuska1989}, by choosing either only odd or only even polynomial orders. $hp$-multigrid methods can be applied on both uniform and refined grids, and generally define the hierarchy of coarse spaces by first performing a $p$-coarsening of the mesh until the $p=1$ level before applying an $h$-coarsening. In the case of uniform meshes, $h$-coarsening is usually applied using a standard $h$-multigrid approach, with the coarsest multigrid level comprising a few large elements e.g.\ \cite{Nastase2006}. Hierarchical $hp$-multigrid approaches are commonly applied on refined grids and often replace the standard $h$-refined linear nodal basis with a hierarchical basis, see e.g.\ \cite{Yserentant1986b,Mitchell2010}. Mitchel et al.\ apply a hierarchical $hp$-approach to refined triangular meshes in \cite{Mitchell2010} and show that the computational complexity of a single V-cycle is $\mathcal{O}(M/p)$, where $M$ is the number of nonzero entries in the stiffness matrix.

The main objective of this contribution is to present a hierarchical multigrid method for the solution of linear systems arising from the finite cell method on both uniform grids and multi-level $hp$-discretizations. The proposed scheme takes advantage of the hierarchical nature of the basis functions in the finite cell method \cite{Parvizian2007,Duster2008} and the multi-level $hp$-scheme \cite{Zander2015,Zander2016}. The FCM implementation used in this work, is based on the $p$-version of the finite element method ($p$-FEM) \cite{Babuska1981} and utilizes high-order integrated Legendre shape functions that are hierarchical in nature, i.e. the set of basis functions of order $p$ contains all basis functions from order 1 up to $p$. Likewise, the superposition principle used to perform spatial refinement in the multi-level $hp$-method results in a hierarchical basis. The suggested multigrid scheme is therefore elegant and efficient, simple in implementation and not strongly interwoven with the code, since all restriction and prolongation operators reduce to binary matrices which do not need to be explicitly constructed. Furthermore, our scheme makes use of additive Schwarz smoothing, and builds upon the work done on additive Schwarz preconditioning for uniform and multi-level $hp$-refined finite cell systems in \cite{Jomo2019}. The overall algorithm results in convergence rates that are independent of the mesh resolution, refinement level and, under certain conditions, the polynomial order. The multigrid technique put forward in this contribution is used to solve large finite cell systems with up to 3.2 billion degrees of freedom (DOFs) and shown to be suitable for large-scale finite cell analyses on massively parallel machines. 

The paper at hand is organized into five sections. The core features of the finite cell method and multi-level $hp$-refinement are highlighted in Section 2. Section 3 presents a hierarchical multigrid framework for FCM and the multi-level $hp$-method focusing on the construction of the coarse spaces and the choice of smoothers to guarantee robustness in the immersed setting. A series of numerical examples are handled in Section 4. These examples are comprised of simple benchmark test cases as well as large-scale applications that attest the suitability of the proposed solution scheme for large-scale finite cell analyses. Section 5 concludes the paper, providing a summary and an outlook.

\renewcommand{\picsDir}{methods/pics}
\renewcommand{\graphDir}{methods/graphs}
\section{The finite cell method and multi-level $hp$-refinement}
The following section presents the fundamental ideas behind the finite cell method and the multi-level $hp$-method. These discretization techniques are used in tandem to perform numerical simulations on domains with complex geometries. Noteworthy application areas include the analysis of bone-implant systems \cite{Elhaddad2017}, the simulation of metal additive manufacturing processes \cite{Ozcan2018} and the modeling of crack growth in brittle structures by means of a phase-field approach \citep{Hug2020}. 

\begin{figure}[t]
    \centering
    \subfloat[Domain of computation.]{ \includegraphics[width=0.33\textwidth]{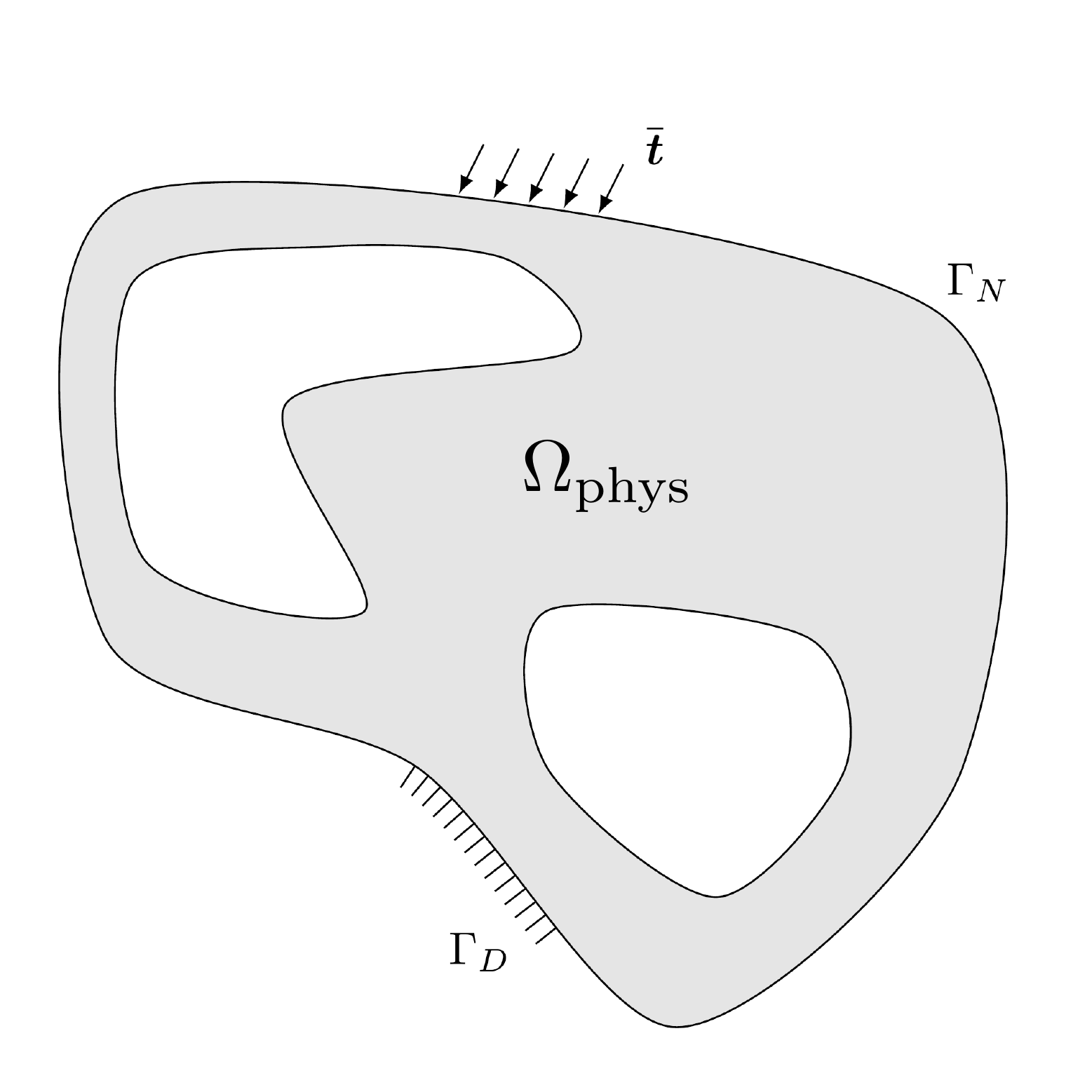} }%
    \subfloat[Fictitious domain extension.]{ \includegraphics[width=0.33\textwidth]{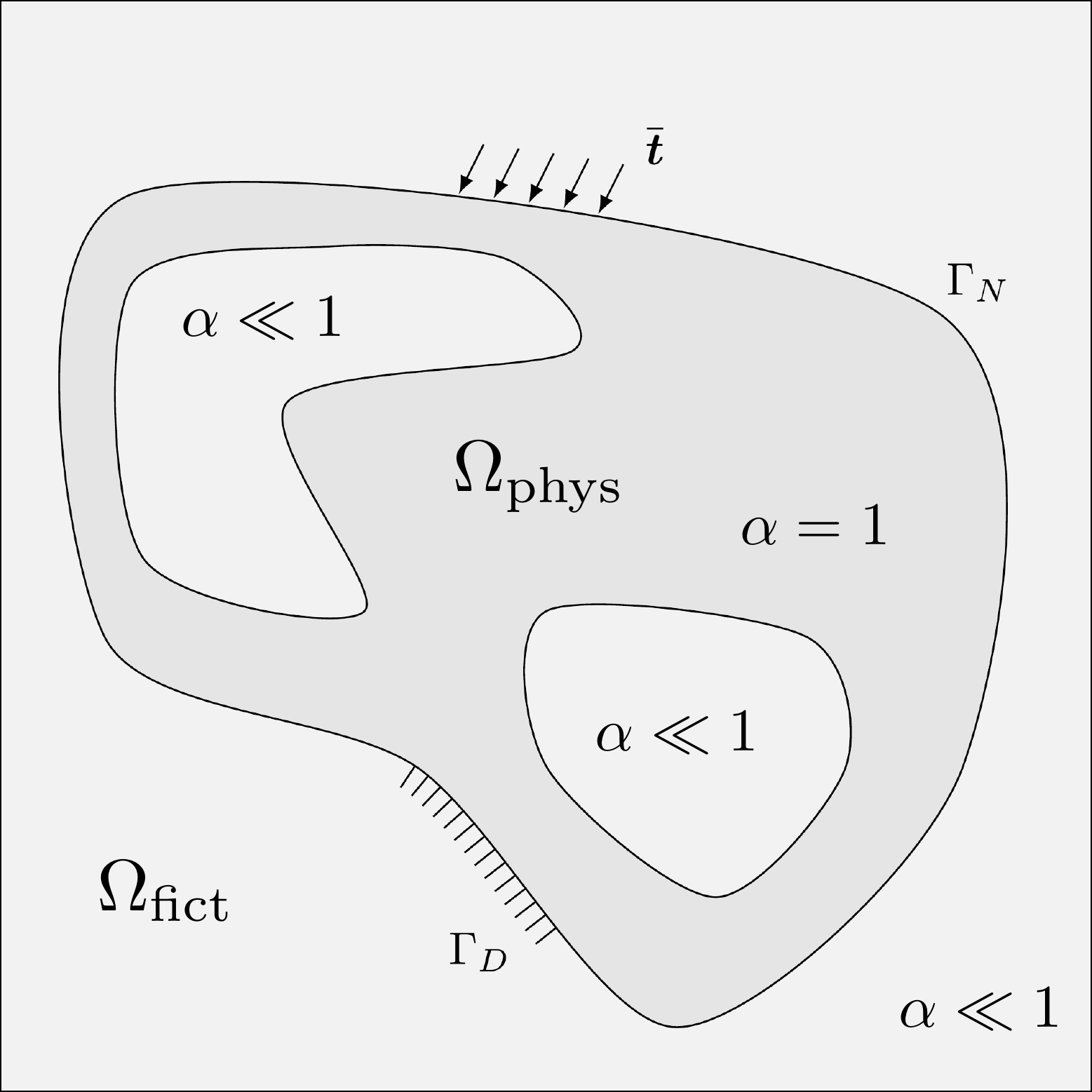} }%
    \subfloat[Structured $hp$-refined FCM mesh.]{ \includegraphics[width=0.33\textwidth]{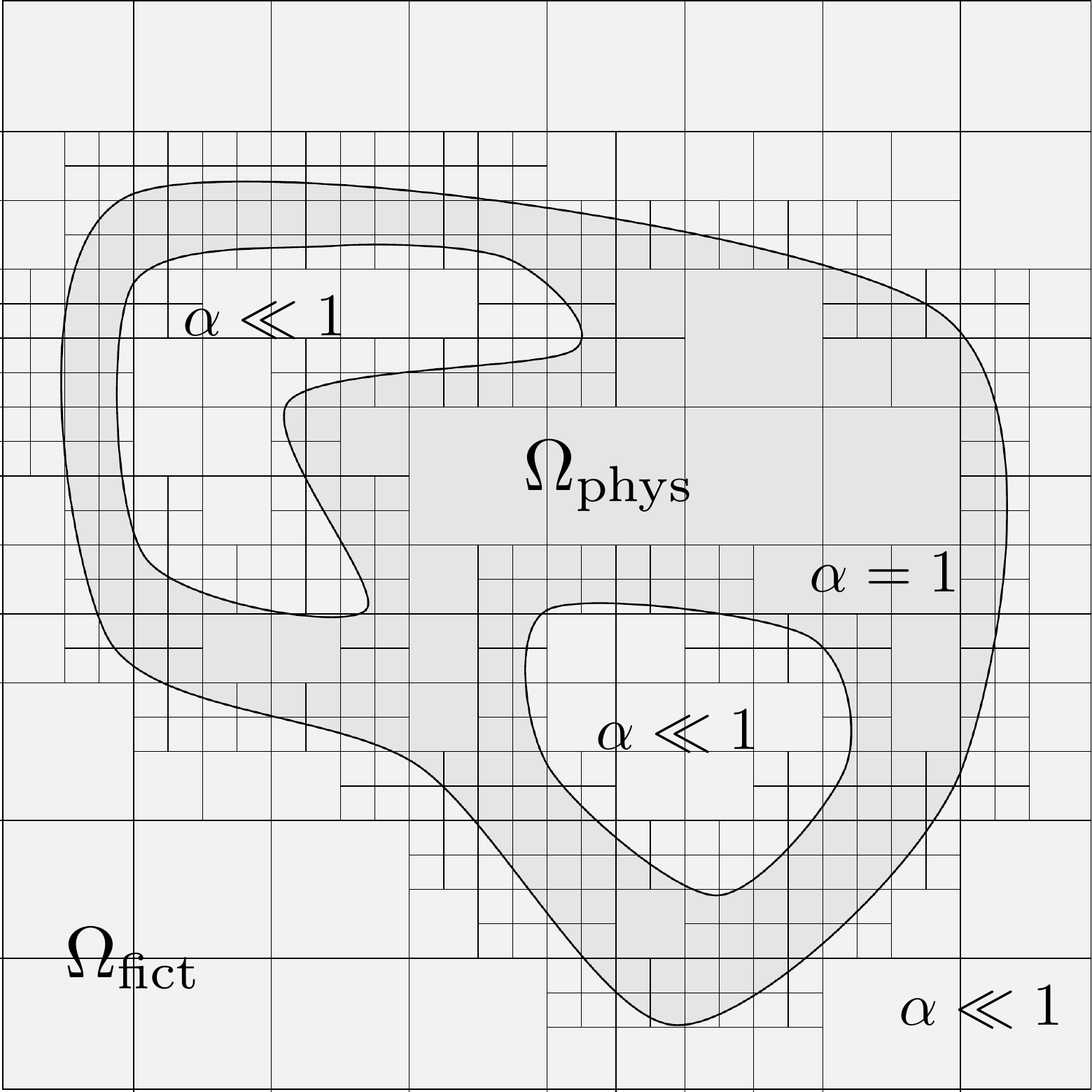} }%
    \caption{Illustration of the main idea behind the finite cell method.}
    \label{fig::fcm}
\end{figure}

\subsection{Ingredients of the finite cell method}
The finite cell method (FCM) is an immersed finite element method that combines a fictitious domain approach with high-order finite elements making it suitable for performing FE-analysis on bodies with very complex geometries. Consider a body defined by the physical domain $\Omega_{\text{phys}}$ and boundary $\Gamma$ with $\Gamma  = \Gamma_D \cup \Gamma_N$ and $\Gamma_D \cap \Gamma_N = \emptyset$. $\Gamma_D$ represents the part of the boundary where Dirichlet boundary conditions $g_D$ are applied while $\Gamma_N$ represents the Neumann boundary with applied traction $g_N$. The main idea of the finite cell method is to place $\Omega_{\text{phys}}$ in an immersing (fictitious) domain $\Omega_{\text{fict}}$ yielding an overall computational domain $\Omega$ of simple shape that can be easily discretized by regular finite elements as illustrated in Figure \ref{fig::fcm}. An indicator function $\alpha$ is commonly used to distinguish between points lying within $\Omega_{\text{phys}}$, which are associated with the value $\alpha = 1$, and points lying in $\Omega_{\text{fict}}$, where alpha is chosen as a small constant $\alpha = \epsilon$ with $\epsilon \ll 1$. This indicator function is applied to the weak form in the FCM as shown in Equation \ref{eq::weakform} and can be interpreted as a penalization of the fictitious domain that results in an asymptotically consistent method recovering the original weak form when $\alpha \rightarrow 0$.
\begin{eqnarray}\label{eq::weakform}
    a(\mathbf{u}_h,\mathbf{v}_h) &=&  \int (\mathbf{B} \mathbf{v}_h )^T \alpha(\boldsymbol{x}) \mathbf{C} \mathbf{B} \mathbf{u}_h \ \textrm{d}\Omega + \int \beta \,(\mathbf{N} \mathbf{v}_h )^T \mathbf{N} \mathbf{u}_h \ \textrm{d}\Gamma_D \nonumber \\
    b(\mathbf{v}_h) &=& \int \alpha(\boldsymbol{x}) (\mathbf{N}\mathbf{v}_h)^T \mathbf{f} \ \textrm{d}\Omega + \int (\mathbf{N} \mathbf{v}_h)^T g_N \ \textrm{d}\Gamma_N + \int \beta \, (\mathbf{N} \mathbf{v}_h)^T \, g_D \ \textrm{d}\Gamma_D 
\end{eqnarray}
Equation \eqref{eq::weakform} shows the discretized weak form of a finite cell problem in the case of linear elasticity with a symmetric bilinear form $a(\mathbf{u}_h,\mathbf{v}_h)$ and linear functional $b(\mathbf{v}_h)$. The term $\mathbf{B}$ denotes the linear strain operator, $\mathbf{N}$ the shape function vector, $\mathbf{C}$ the elasticity tensor and $\boldsymbol{f}$ a prescribed volumetric force. The unknown displacement field vector and test functions are represented by $\mathbf{u}_h$ and $\mathbf{v}_h$, respectively. The penalty method \cite{Babuska1973} is used in \eqref{eq::weakform} to impose the Dirichlet boundary conditions and $\beta$ denotes the penalty parameter.

Various kinds of embedding meshes can be employed in the finite cell method. We use a mesh based on $p$-FEM that utilizes integrated Legendre polynomials. The hierarchical nature of these basis functions can be exploited to construct a multigrid algorithm as outlined in Section \ref{sec::multigrid}.

The fictitious domain approach employed in the finite cell method has two main benefits. First, it significantly simplifies the meshing process, as the generation of (high-order) body-fitting meshes can be especially challenging for complex geometries in three dimensions. Secondly, it yields a versatile discretizational framework that can be applied in a straightforward manner to different types of geometric models such as those stemming from CAD systems, constructive solid geometry, scan images and even oriented point clouds. The fictitious extension, however, has several implications that make the simulation pipeline in FCM, and immersed methods in general, differ from that of conventional body-fitted methods. Three core aspects that require special treatment in immersed methods include \textit{i)} the numerical integration on cut elements, \textit{ii)} the application of Dirichlet boundary conditions and \textit{iii)} the conditioning of the system. A comprehensive description of each of these areas is beyond the scope of this article, as it focuses only on the conditioning and solution of the resulting linear systems. The interested reader is directed to the review articles on the finite cell method \cite{Schillinger2014,Duster2017} for a comprehensive overview of cut cell integration techniques and the weak imposition of Dirichlet boundary conditions.

\subsection{The multi-level $hp$-method}

\begin{figure}[t]
  \begin{minipage}[c]{0.7\textwidth}
    \begin{minipage}[t]{0.45\textheight}	  
     \begin{center}
       \includegraphics[width=0.95\textwidth]{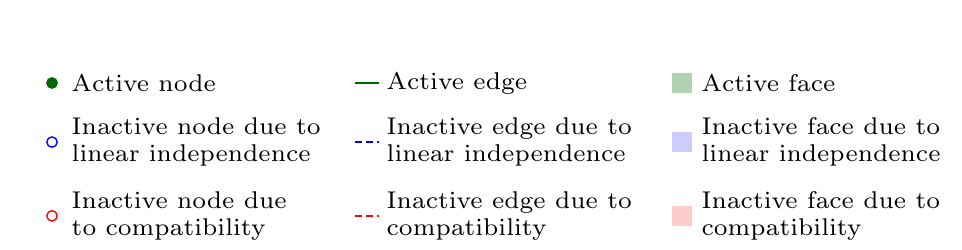}
     \end{center}	
    \end{minipage}	  
    \begin{minipage}[t]{0.5\textheight}
     \centering	    
     \begin{center}
       \subfloat[One-dimensional case]{ \includegraphics[width=0.56\textwidth]{pics/refinement_1d_mlhp_woNumbers} }
       \subfloat[Two-dimensional case]{ \includegraphics[width=0.38\textwidth]{pics/multiLevelhp2d.tikz} }%
     \end{center}	
   \end{minipage}	  
  \end{minipage}
  \begin{minipage}[c]{0.28\textwidth}
   \centering	 
   \vspace{4mm} 
   \begin{center}
      \includegraphics[width=0.95\textwidth]{pics/multiLevelhp3d-notext.tikz}%
    \end{center}
    \footnotesize
    (c) Three-dimensional case
  \end{minipage}
  \caption{Illustration of the multi-level $hp$-refinement scheme with two refinement levels, $k=2$, in different spatial dimensions. The deactivation of specific topological components following a simple rule-set ensures compatibility and linear independence of the basis functions  \cite{Zander2016}.}
  \label{fig::multilevelhp::multilevelhpIdea}%
\end{figure}

The multi-level $hp$-method presented in \cite{Zander2015,Zander2016} is a novel $hp$-method that performs spatial refinement based on superposition. The scheme was
developed to circumvent the challenges associated with constraining arbitrary levels of hanging nodes, while at the same time preserving the desirable characteristics of classical $hp$-formulations. It employs high-order overlay meshes which are placed over coarse elements in critical areas of the mesh. The original formulation of the $hp$-scheme is  based on the $p$-version of the finite element method and makes use of the direct association of the topological components --- nodes, edges, faces and solids --- with the degrees of freedom to ensure the compatibility and linear independence of the basis functions. The routines used to construct the basis are intuitive and make the method suitable for dealing with complex refinement scenarios in multiple dimensions. This contribution only considers refinement patterns that are driven by a priori geometrical information. The use of error estimators to guide refinement in the scheme is treated in \cite{DAngella2016}. A recent publication presents a novel $hp$-adaptive strategy for multi-level discretiztations \cite{Darrigrand2020}. Figure \ref{fig::multilevelhp::multilevelhpIdea} illustrates the manner in which topological components can be activated and deactivated in order to construct the multi-level $hp$-basis. 

A multi-level $hp$-mesh comprises elements with different refinement levels as shown in Figure \ref{fig::multilevelhp::multilevelhpIdea}. The letter $k$ is used to denote the refinement level/depth of an element. Elements on the lowest refinement level with $k=0$ are termed \textit{base elements}. Spatial refinement is performed by superposing a base element with subelements formed by uniform bisections. This procedure is performed recursively until the desired refinement depth, resulting in a refinement tree for every refined base element. The term \textit{leaf elements} is used to refer to the set of all elements in the mesh that do not have subelements. An in-depth description of the properties of the multi-level $hp$-basis and the steps needed for its construction can be found in \cite{Zander2016,Zander2016b}.

\renewcommand{\picsDir}{multigrid/pics}
\renewcommand{\graphDir}{multigrid/graphs}
\section{Hierarchical multigrid for uniform and multi-level $hp$-grids}\label{sec::multigrid} 

The following sections handle the use of multigrid algorithms for the solution of finite cell systems on uniform meshes and multi-level $hp$-refined grids. It begins with a brief summary of the multigrid method, introducing the notation and terminology that is employed in this contribution. It thereafter expounds on how an $hp$-multigrid method can be devised for finite cell systems involving multi-level $hp$-refinement and elaborates on suitable smoothing techniques for achieving convergence rates that are independent of the cut configuration, mesh size and in special cases, even the polynomial order.  

\subsection{The multigrid algorithm}
Multigrid methods aim at improving the convergence of a linear system by incorporating information obtained from a hierarchy of coarse discretizations. Consider a sequence of function spaces $\mathcal{V}_0, \mathcal{V}_1 \dots \mathcal{V}_{\ell-1}, \mathcal{V}_{\ell}$, where the index $\ell$ denotes the level-number and $\ell \in [0, \ell_{max}]$. $\mathcal{V}_{\ell_{max}}$ denotes the finest space which contains the solution of the original problem, whereas $\mathcal{V}_0$ denotes the coarsest space. Each multigrid level is associated with a system of linear equations  
\begin{equation}
    \mathbf{A}_{\ell} \mathbf{x}_{\ell} = \mathbf{b}_{\ell}.
\end{equation}
In this work, only sequences of nested spaces are considered, with the consequence that basis functions in a coarse space can be expressed as a linear combination of basis functions from a finer space. This makes it possible to define restriction operators $\mathbf{R}_{\ell}$ that map basis functions from the fine space of level $\ell$ into the next coarse space of level $\ell - 1$. Conversely, the prolongation operators $\mathbf{R}_{\ell}^T$ can be introduced that map a vector in the coarse space $\ell -1$ to a vector in the fine space $\ell$ which corresponds to the same function. These operators allow quantities such as vectors and matrices to be transferred from one level to another e.g.\ $\mathbf{A}_{\ell-1} = \mathbf{R}_{\ell}\mathbf{A}_{\ell}\mathbf{R}_{\ell}^T$ and $\mathbf{r}_{\ell-1} = \mathbf{R}_{\ell}\mathbf{r}_{\ell}$, where $r_{\ell}$ and $r_{\ell -1}$ are the residuals on level $\ell$ and $\ell -1$, respectively.

The two main steps in a multigrid iteration are a smoothing process and the coarse-grid correction. These steps are considered complementary since they act on different components of the error. The error in each multigrid level can be expressed as the difference between the exact solution $\mathbf{x}_{\ell}$ and the current approximation $\tilde{\mathbf{x}}_{\ell}$, i.e.\
\begin{equation}
    \mathbf{e}_{\ell} = \mathbf{x}_{\ell} - \tilde{\mathbf{x}}_{\ell}
\end{equation}
The smoothing process, or smoothing in short, efficiently reduces the oscillatory components of the error that are associated with large eigenvalues in level $\ell$. Smoothing is performed by the repeated application of linear iterative methods such as fixed-point iterations. Standard smoothers include the weighted Jacobi method and the Gauss-Seidel method, but it is also possible to use fixed-point smoothing based on additive Schwarz (AS) techniques, multiplicative Schwarz techniques and incomplete LU factorizations. In this work, the smoothing process is based on additive Schwarz techniques since they can be readily applied in massively parallel settings and their construction can be done in a fully parallel manner without the need for synchronization. This is in contrast to techniques such as Gauss-Seidel and multiplicative Schwarz, which often require coordinated application in parallel settings such as multicoloring algorithms, see e.g.\ \citep{Saad2003}. The application of the smoother is carried out per the formula
\begin{equation}
    \tilde{\mathbf{x}}_{\ell} = \tilde{\mathbf{x}}_{\ell} + \omega \mathbf{M}_{\ell}^{-1} \mathbf{r}_{\ell},
\end{equation}
where $\omega$ denotes the relaxation parameter and $\mathbf{M}^{-1}$ the smoother. The choice of the relaxation parameter $\omega$ for the various FCM meshes applied in this manuscript is elaborated in Section \ref{sec::smoothingmlhpfcm}. The coarse-grid correction accelerates the convergence of the smooth components of the error that are associated with small eigenvalues in level $\ell$. By taking advantage of the relation between the residual $\mathbf{r}_{\ell}$ and the error $\mathbf{e}_{\ell}$ 
\begin{equation}
    \mathbf{A}_{\ell}\mathbf{e}_{\ell} = \mathbf{r}_{\ell}
\end{equation}
it is possible to approximate the smooth components of the error $\mathbf{e}_\ell$ using a coarser discretization (grid), where $\mathbf{e}_{\ell} \approx \mathbf{R}_{\ell}^T \mathbf{e}_{\ell - 1}, $ by solving the coarse-grid system
\begin{equation}
    \mathbf{A}_{\ell-1}\mathbf{e}_{\ell-1} = \mathbf{r}_{\ell-1}.
\end{equation}
The righthand side vector of this coarse system is formed by restricting the residual at level $\ell$ to level $\ell -1$, i.e.\ $\mathbf{r}_{\ell - 1} = \mathbf{R}_{\ell} \mathbf{r}_{\ell}$. Once this system has been solved, its solution, termed the \textit{coarse-grid correction}, can be transferred to level $\ell$ through the expression 
\begin{equation}
    \tilde{\mathbf{x}}_{\ell} = \tilde{\mathbf{x}}_{\ell} + \mathbf{R}_{\ell}^T\mathbf{e}_{\ell-1}.
\end{equation}

\begin{figure}[H]
    \centering
  \includegraphics[width=0.99\textwidth]{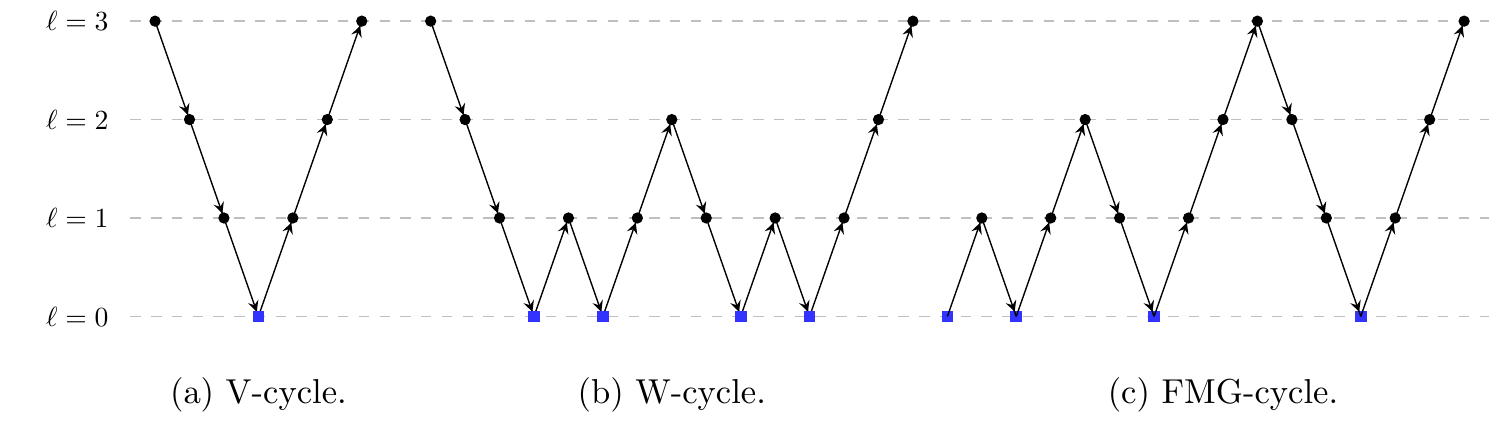}
    \caption{Illustration of three different types of multigrid iterations. The circular dots represent $n_s$ smoothing steps performed on every level with $\ell \neq 0$, while the square dots represent an \enquote{exact} solve performed on $\ell = 0$. Downwards pointing arrows represent restriction operations, while upwards pointing arrows represent prolongations.  }
\label{fig::multigridcycles}
\end{figure}

A single multigrid iteration or cycle consists of the application of the smoothing and coarse-grid correction steps on the different multigrid levels. An exception is, however, made on the coarsest level $\ell = 0$. Here, the coarse system is solved with either a direct solver, if the system is small in size, or with an iterative solver before prolongating its solution to level $\ell = 1$. There are different ways of traversing the levels within a multigrid cycle, the most prevalent of which are the V-cycle, W-cycle and the full multigrid cycle (FMG). A schematic representation of each of these cycles is given in Figure \ref{fig::multigridcycles}. V-cycles are applied in this work and a summary of the multigrid algorithm for this iteration type is given in Algorithm \ref{alg::MGVcycle}. 

\begin{algorithm}
    \caption{ $\tilde{\mathbf{x}}_{\ell}$ = performVcycle($\tilde{\mathbf{x}}_{\ell},\mathbf{r}_{\ell},{\ell}$)}
    \label{alg::MGVcycle}
     \vspace{.2cm}
    \eIf { $l \neq 0$ }
    {
     \textcolor{blue}{\# perform $n_s$ pre-smoothing steps } \\
     \For{$i \in n_s$ }
     {
         $\tilde{\mathbf{x}}_{\ell} \leftarrow \tilde{\mathbf{x}}_{\ell} + \omega \mathbf{M}^{-1}_{\ell} \mathbf{r}_{\ell}$ \\
     }
     \vspace{.2cm}
     \textcolor{blue}{\# update residual} \\
     $\mathbf{r}_{\ell} = \mathbf{b}_{\ell} - \mathbf{A}_{\ell} \tilde{\mathbf{x}}_{\ell}$ \\
     \vspace{.2cm}
     \textcolor{blue}{\# coarse grid correction} \\
     $\mathbf{r}_{{\ell}-1} = \mathbf{R}_{\ell} \mathbf{r}_{\ell}$ \\
     $\tilde{\mathbf{x}}_{{\ell}-1} = \text{performVcycle}(\mathbf{0}, \mathbf{r}_{{\ell}-1}, {\ell} - 1 )$ \\
     $\tilde{\mathbf{x}}_{{\ell}} = \tilde{\mathbf{x}}_{{\ell}} + \mathbf{R}^{T}_{\ell} \tilde{\mathbf{x}}_{{\ell}-1} $ \\
     $\mathbf{r}_{\ell} = \mathbf{b}_{\ell} - \mathbf{A}_{\ell} \tilde{\mathbf{x}}_{\ell} $ \\

     \vspace{.2cm}
     \textcolor{blue}{\# perform $n_s$ post-smoothing steps } \\
     \For{$i \in n_s$ }
     {
         $\tilde{\mathbf{x}}_{\ell} \leftarrow \tilde{\mathbf{x}}_{\ell} + \omega \mathbf{M}^{-1}_{\ell} \mathbf{r}_{\ell}$ \\
     }
    }
    {
     \textcolor{blue}{\# solve the coarse system} \\
     $\tilde{\mathbf{x}}_{\ell} = \text{solve}(\mathbf{A}_{\ell},\mathbf{r}_{\ell})$
    }
     \vspace{.2cm}
\end{algorithm}

\subsection{Multigrid for hierarchical bases}
It is possible to take advantage of the hierarchical nature of the shape functions in the finite cell method and multi-level $hp$-refinement to devise a hierarchical multigrid method. 
Let $\mathbf{N}_{\ell}$ denote the basis functions that belong to the multigrid level $\ell$. For each level in a $p$-multigrid (for uniform meshes) or $hp$-multigrid (for multi-level $hp$-meshes) the following relation holds  
\begin{equation}\label{eq::hierarchicalFunctions}
    \mathbf{N}_0 \subset \mathbf{N}_1 \dots \mathbf{N}_{\ell - 1} \subset \mathbf{N}_{\ell}.
\end{equation}
The hierarchical structure in \eqref{eq::hierarchicalFunctions} is also reflected in the degree of freedom vector $\mathbf{u}$. The vector of DOFs on level $\ell$, denoted by $\mathbf{u}_{\ell}$, is made up of coefficients from level $\ell - 1$, represented by $\mathbf{u}_{\ell - 1}$ and entries solely on level $\ell$, denoted by $w_{\ell}$, i.e.\
\begin{equation}
    \mathbf{u}_{\ell} = \left[ 
\begin{array}{c}
    \mathbf{u}_{\ell - 1} \\
    \mathbf{w}_{\ell}
\end{array} \right].
\end{equation}
From \eqref{eq::hierarchicalFunctions} it follows that the basis functions spanning the (coarse) nested subspace can not only be constructed by a linear combination of basis functions in the fine space, which is done in restriction and prolongation, but are explicitly contained in the basis functions of the fine space. This leads to an elegant and efficient multigrid framework since all restriction and prolongation operators reduce to binary matrices which do not need to be explicitly applied. Transitioning from one multigrid level to another is done easily by either leaving out specific basis functions or re-introducing them back into the system. Equation \eqref{eq::restrictionoperation} illustrates how the DOFs on level $\ell$ can be \enquote{trimmed} to obtain the DOFs on level $\ell -1$. $\mathbf{I}$ represents an identity matrix while the term $\mathbf{0}$ denotes a matrix in which all entries are zero.
\begin{equation}\label{eq::restrictionoperation}
    \mathbf{u}_{\ell - 1} = \left[ \mathbf{I}, \ \mathbf{0} \right] \left[ 
\begin{array}{c}
    \mathbf{u}_{\ell - 1} \\
    \mathbf{w}_{\ell}
\end{array} \right].
\end{equation}
For uniform grids with high-order elements, an \textit{arithmetic $p$-sequence} is used to generate the coarse subspaces, i.e.\ the value of $p$ is progressively reduced until $p=1$, see \cite{Mitchell2010}. In the case of multi-level $hp$-grids, it is required to first reduce the polynomial order $p$ in the overlay meshes, before reducing the levels of refinement. This procedure ensures that the resulting coarse spaces are subspaces.

\begin{figure}[H]
    \centering
  \includegraphics[width=0.95\textwidth]{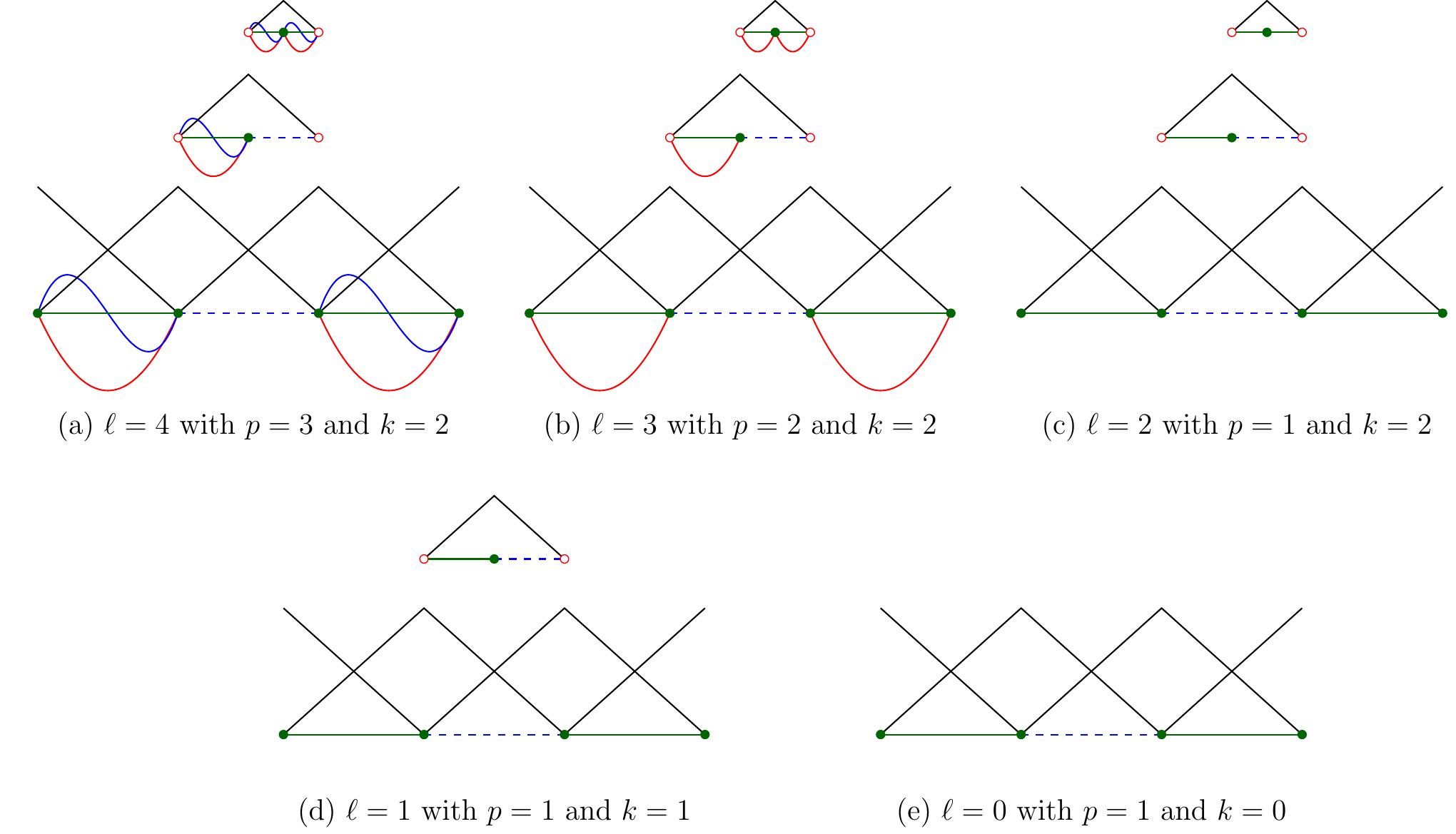}
    \caption{Multigrid levels in a one-dimensional mesh with two levels of multi-level $hp$-refinement ($k=2$) and a polynomial order of $p=3$.}
    \label{fig::hpmultigrid}
\end{figure}

The hierarchical nature of the FE-basis in FCM and the multi-level $hp$-scheme is also reflected in the system matrices. Consequently, computation of lower-level matrices is not needed as this information is readily available. Equation \eqref{eq::hierarchicalMatrices} shows the structure of a hierarchical matrix of level $\ell$, which consists of entries $\tilde{\mathbf{A}}_{\ell}$ belonging to basis functions contained solely in the highest level, a term $\tilde{\mathbf{A}}_{\ell-1}$ that contains entries of all lower levels up to $\ell -1$ and a term $\tilde{\mathbf{A}}_{\ell,\ell-1}$ that couples DOFs on level $\ell$ with all other DOFs. 
\begin{equation}\label{eq::hierarchicalMatrices}
    \mathbf{A}_{\ell} =\left[
\begin{array}{cc}
    \tilde{\mathbf{A}}_{\ell} & \tilde{\mathbf{A}}_{\ell,\ell-1} \\
    \tilde{\mathbf{A}}_{\ell,\ell-1}^T & \mathbf{A}_{\ell-1} \\
\end{array} \right]
\end{equation}
No distinction will be made from this point on in the manuscript between the $p$-multigrid scheme for uniform meshes and the $hp$-multigrid approach of the multi-level $hp$-discretizations, since the $p$-multigrid can be regarded as a special case of the $hp$-multigrid in which $k=0$. Figure \ref{fig::hpmultigrid} shows the $hp$-multigrid levels used for a one-dimensional multi-level $hp$-mesh. A reduction in the $p$-level is performed first, followed by a reduction in the levels of refinement until the lowest multigrid level with $p=1$ and $k=0$ is reached.

\begin{remark}
    \textit{The $p=1$ and $k=0$ level is chosen as the coarsest level in the proposed $hp$-multigrid approach. This allows a simple yet elegant parallel implementation, since only the DOFs on the finest problem need to be distributed over the MPI tasks. All coarse sub-problems in the multigrid hierarchy are able to reuse the parallel data structures that have already been set up for the fine problem and do not require any additional steps such as redistribution or partitioning of the DOFs. It should be noted, however, that the $p=1$, $k=0$ problem can be coarsened further using a standard geometric multigrid algorithm. This procedure can improve the convergence of the coarse problem and is beneficial in applications where the solution of the coarse problem is slow. }
\end{remark}

\subsection{Selection of suitable smoothing strategies}\label{sec::smoothingmlhpfcm}
The choice of an appropriate smoother in a multigrid algorithm is fundamental in achieving convergence rates that are independent of the mesh parameter $h$. In specific problems, it is even possible to achieve convergence rates that are independent of the element polynomial order $p$ when a suitable smoother is chosen, see Section \ref{sec::poissonmeshfitting}. 

\subsubsection{Additive Schwarz smoothing}
\begin{figure}[H]
    \centering
    \subfloat[Smallest eigenmode of a boundary-conforming system with Jacobi preconditioning.]
    {
      \includegraphics[width=0.40\textwidth]{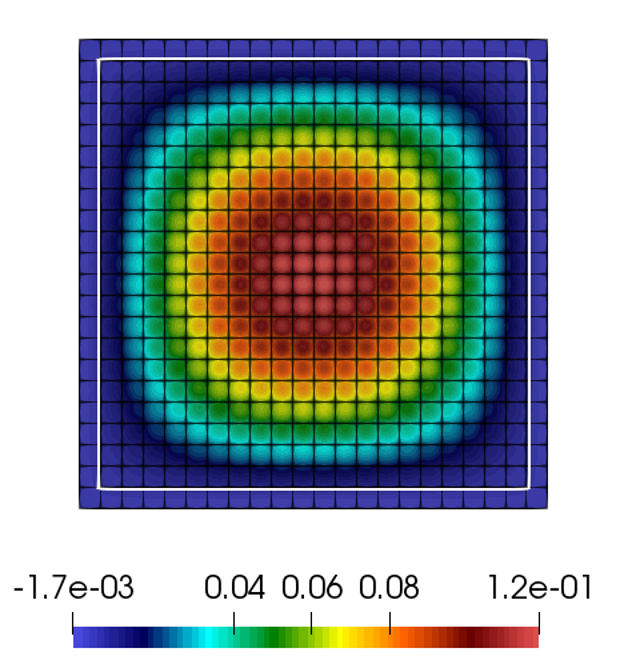}
      \label{fig::smJacobiBC}
    }
    \hspace{5mm}
    \subfloat[Smallest eigenmode of an immersed system with Jacobi preconditioning.]
    {
      \includegraphics[width=0.47\textwidth]{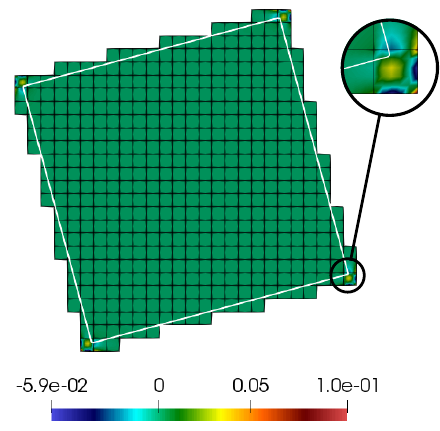}
      \label{fig::smJacobiImmersed}
    }
    \caption{Illustration of the smallest eigenmodes of a Jacobi preconditioned system arising from the Poisson equation on a boundary-conforming grid and an immersed grid.}
    \label{fig::smallesteigenvalues}
\end{figure}

Smoothing techniques such as the Jacobi and Gauss-Seidel methods are commonly used in multigrid algorithms for boundary-conforming finite element methods. These techniques are, however, not suitable for FCM problems as they fail to resolve the conditioning problems associated with cut cells. To illustrate this point, we consider a system arising from the Poisson equation posed on a square domain and compare two scenarios; a case in which the domain is discretized in a boundary-conforming manner and a scenario where an immersed grid is utilized. Figure \ref{fig::smallesteigenvalues} visualizes the eigenmodes associated with the smallest eigenvalue (further denoted as \enquote{smallest eigenmode}) of the system matrix when Jacobi preconditioning is applied. In the boundary-conforming case, the smallest eigenmode corresponds to a smooth mode that spans the entire domain as shown in Figure \ref{fig::smJacobiBC}. This mode can be sufficiently approximated on a coarse space. In the immersed case, the smallest eigenmode is restricted to a few cut elements, see Figure \ref{fig::smJacobiImmersed}. The occurrence of such small modes in immersed methods impedes the convergence of iterative solvers and makes conventional smoothers such as the Jacobi and Gauss-Seidel methods unsuitable for immersed grids. 

Reference \citep{Prenter2019b} shows that smoothers based on the additive and multiplicative Schwarz lemmas are better suited for multigrid methods involving cut cells and presents results of multiplicative Schwarz smoothing in different linear elastic examples. As mentioned in the previous section, an additive Schwarz approach is used in this work, since it can be easily applied in parallel and has been successfully used for preconditioning FCM systems involving multi-level $hp$-refinement, see \cite{Jomo2019}. The core idea behind additive Schwarz smoothing is to group the basis functions in a mesh into \textit{blocks/groups} and thereafter construct a smoother by inverting and summing sub-matrices of $\mathbf{A}$ devised from the chosen blocks as per the expression
\begin{equation}\label{eq:AdditiveSchwarz}
    \mathbf{M}^{-1} = \sum_{i=1}^{n_{\text{blocks}}} \mathbf{P}_i \underbrace{\left(\mathbf{P}_i^T \mathbf{A}\mathbf{P}_i\right)^{-1}}_{\mathbf{A}_i^{-1}}\mathbf{P}_i^T \, ,
\end{equation}
where $n_{\text{blocks}}$ denotes the number of blocks, $i$ the index corresponding to the $i^{\text{th}}$ block containing $m$ basis functions such that $\mathbf{A}_i \in \mathbb{R}^{m \times m}$. $\mathbf{P}$ and $\mathbf{P}^T$ are restriction and prolongation operators and $\mathbf{P}_i \in \mathbb{R}^{n_{\text{DOFs}} \times m }$.

\begin{figure}[H]
    \centering
    \subfloat[Elementwise blocks]
    {
      \label{fig::elementBlocks}
      \includegraphics[width=0.30\textwidth]{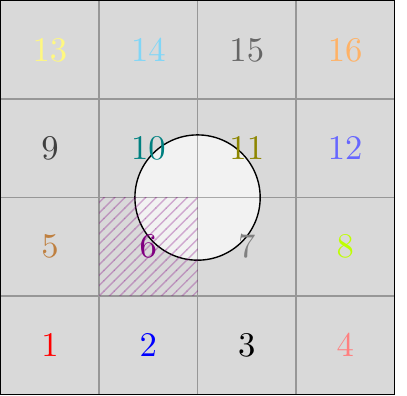}
    }
    \hspace{5mm}
    \subfloat[Patchwise blocks]
    {\label{fig::patchBlocks}
      \includegraphics[width=0.30\textwidth]{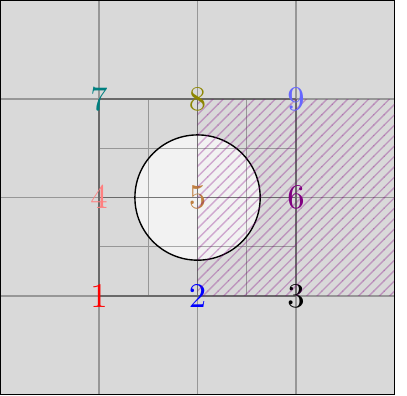}
    }
    \caption{Selection of the additive Schwarz groups for finite cell problems on uniform grids and multi-level $hp$-refined grids. This selection is performed on the granularity of the base elements for uniform grids and node-based element patches for multi-level $hp$-grids. The numbers in the figures indicate the number of additive Schwarz blocks and the support of the $6^{\textrm{th}}$ block is shaded}
\end{figure}

The proper selection of the basis function blocks is essential for guaranteeing a robust smoother that can deal with conditioning problems due to cut elements. Following \cite{Prenter2019,Jomo2019} we select the additive Schwarz blocks in such a way that basis functions that can become potentially linear dependent are grouped together. An \textit{elementwise} approach is used for block selection in uniform grids, i.e.\ all basis functions that are supported on an element are considered to be one group, see Figure \ref{fig::elementBlocks}. A \textit{patchwise} approach is employed for multi-level $hp$-grids, which forms additive Schwarz blocks comprising all basis functions supported on a group of base elements around a mesh node as indicated in Figure \ref{fig::patchBlocks}. It should be noted that these two approaches result in overlapping additive Schwarz blocks. 

\begin{figure}[H]
    \centering
    \subfloat[Smallest eigenmode of an immersed system on a uniform grid with elementwise additive Schwarz preconditioning.]
    {
      \label{fig::imuni}
      \includegraphics[width=0.40\textwidth]{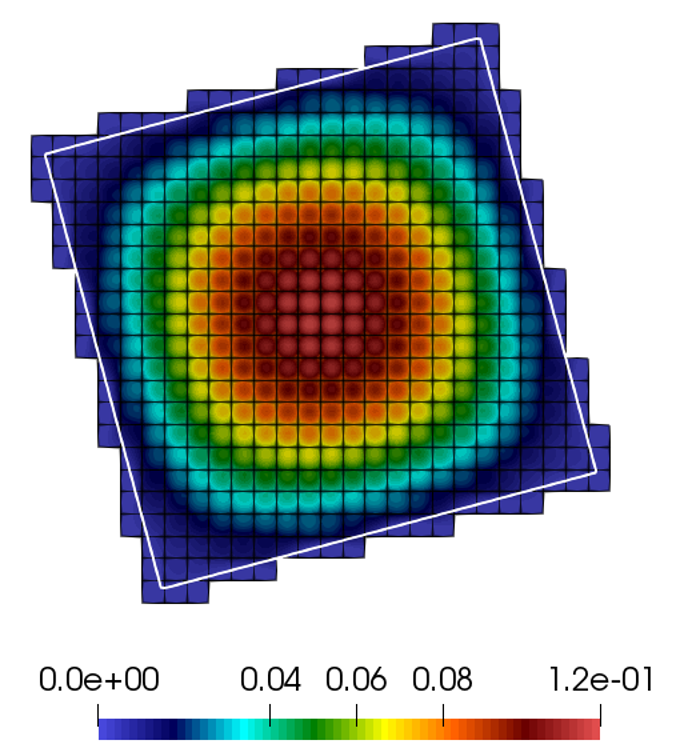}
    }
    \hspace{5mm}
    \subfloat[Smallest eigenmode of an immersed system on a multi-level $hp$-grid with patchwise additive Schwarz preconditioning.]
    {
      \label{fig::immlhp}
      \includegraphics[width=0.39\textwidth]{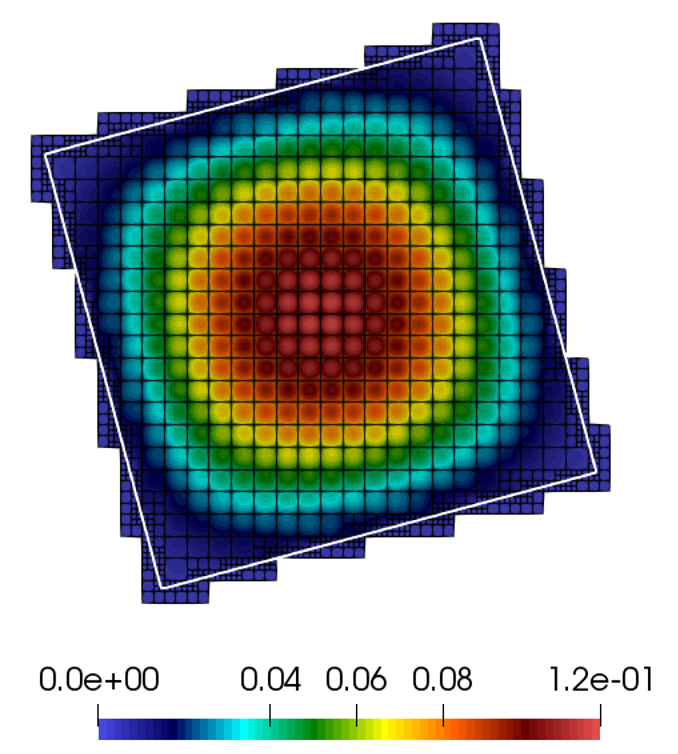}
    }
    \caption{Illustration of the smallest eigenmodes of a system arising from the Poisson equation on a uniform grid, Figure \ref{fig::imuni}, and a multi-level $hp$-grid, Figure \ref{fig::immlhp}. The system matrices are preconditioned using the additive Schwarz techniques. }
    \label{fig::rotatedSmallesteigenvalues}
\end{figure}

We again illustrate the suitability of the proposed smoothers by considering the smallest eigenmodes of a Poisson problem on a square domain that is rotated with respect to a fixed background grid. Figure \ref{fig::rotatedSmallesteigenvalues} visualizes the smallest eigenmodes of the systems that are preconditioned using the elementwise and patchwise additive Schwarz techniques on uniform grids and multi-level $hp$-refined grids, respectively. Both approaches are able to adequately deal with small modes due to cut cells  

\subsubsection{The relaxation parameter $\omega$}

Although the additive Schwarz smoother can be easily applied in a massively parallel setting, it requires sufficient stabilization in order to converge. It is well known that the convergence of fixed-point iterations requires the eigenvalues of the iteration matrix to be bounded, i.e.
\begin{equation}
    \rho( \mathbf{I} - \omega \mathbf{M}^{-1}\mathbf{A} ) < 1.
\end{equation}
In \citep{Prenter2019b}, it is shown that the largest eigenvalue of the matrix $\mathbf{M}^{-1}\mathbf{A}$ is bounded from above by the maximum overlap $n_{\text{max}}$ of the additive Schwarz blocks i.e.\ $\lambda_{\text{max}}(\mathbf{M}^{-1}\mathbf{A}) \leq n_{\text{max}}$. Since the value of $n_{max}$ is $2^d$ for Cartesian grids, the relaxation parameter $\omega$ can be chosen as $\omega = 2/n_{\text{max}}$, thus guaranteeing stability of the fixed-point iteration. Note, however, that the relation between the largest eigenvalue and $n_{\text{max}}$ is an inequality and that it is possible to choose values for $\omega$ that are higher than $2/n_{\text{max}}$. For the integrated Legendre basis functions considered in this work, the best performance was achieved by choosing $\omega=1/3$ for two-dimensional problems and $\omega = 2/15$ for three-dimensional problems, when using an additive Schwarz smoother based on elementwise blocks. When the additive Schwarz blocks are chosen in a patchwise manner, $n_{\text{max}}$ is larger, and the best convergence rates are achieved using $\omega=1/6$ and $\omega=1/18$ in two and three dimensions respectively. The values of $\omega$ suggested in this work were determined in a heuristic approach and achieve the mesh-independent convergence rates for different geometries as shown in Section \ref{sec::numerics}.

It should be noted that different smoothing strategies can be applied on the different multigrid levels, e.g.\ additive Schwarz on the finest level and Gauss-Seidel or Jacobi smoothing on lower levels. This is, however, not investigated in this work and additive Schwarz smoothing is applied on all levels with $\ell \neq 0$.

\subsubsection{Implementational aspects}
\begin{figure}[H]
    \centering
    \subfloat[Maximum AS group size for uniform grids.]
    {
      \label{fig::mUniform}
      \includegraphics[width=0.48\textwidth]{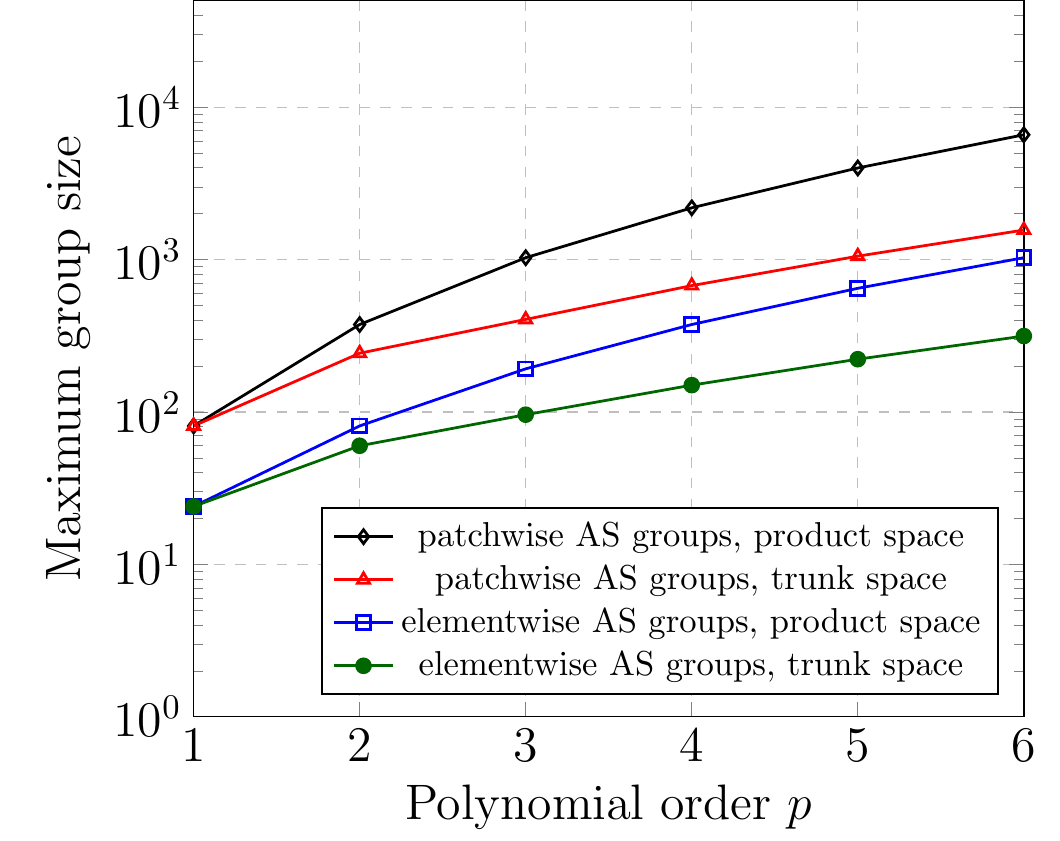}
    }
    \hfill
    \subfloat[Maximum group size for the benchmark example involving multi-level $hp$-refinement in Section \ref{sec::cubeSpheresEx}.]
    {
      \label{fig::refUniform}
      \includegraphics[width=0.48\textwidth]{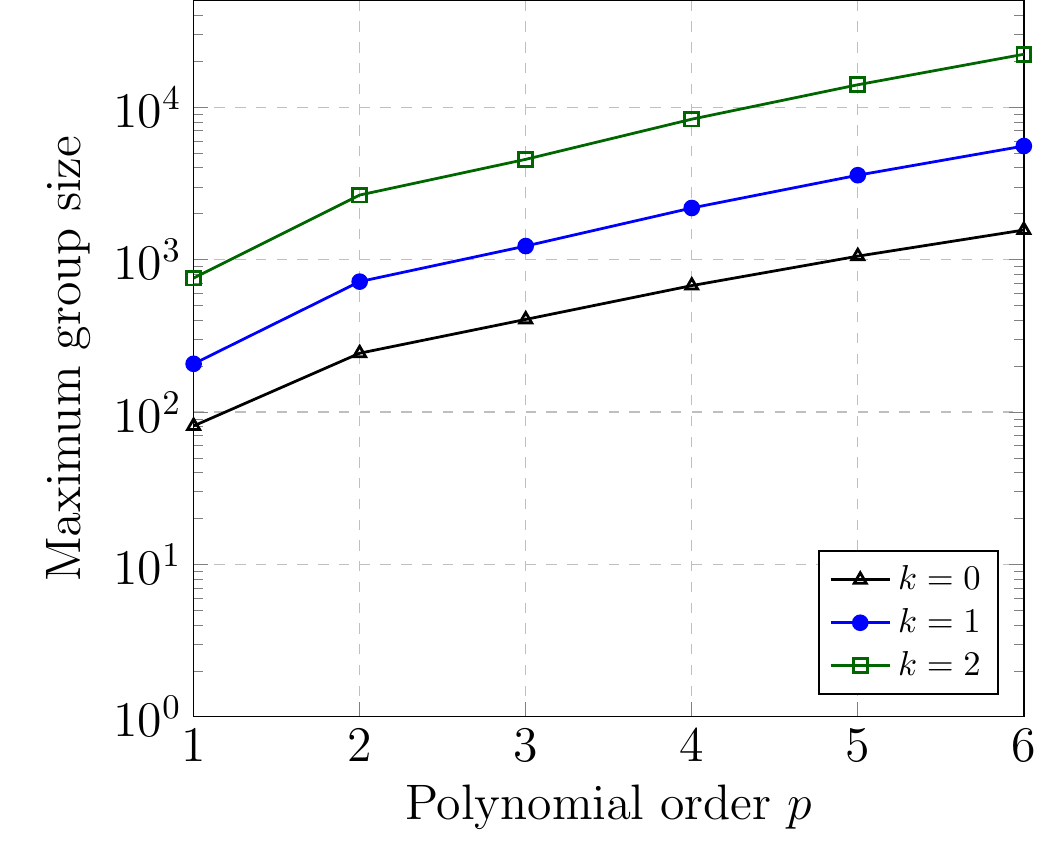}
    }
    \caption{Illustration of the maximum size of the basis function groups used to construct the additive Schwarz smoothers for three-dimensional linear elastic problems.}
    \label{fig::boundsGraph}
\end{figure}

    The multigrid scheme proposed in this paper is implemented in an in-house finite element code written in \texttt{C++}. The code's modular design enables it to be used on a variety of computing platforms ranging from desktop computers to massively parallel systems. It supports hybrid parallel simulations through the use of MPI \cite{mpi1994} and OpenMP\cite{openmp08}. Generation of the computational mesh in large-scale computations is carried out in a distributed manner using a parallel adaptive Cartesian grid. This strategy ensures that the memory resources are utilized in a scalable manner, since no single MPI task needs to know the complete extent of the computational domain. The code framework utilizes the functionality of the \texttt{Epetra} package in \texttt{Trilinos} \cite{Trilinos} for the parallel construction and storage of the system matrices. \texttt{Epetra} is also used to perform parallel linear algebra operations. The system matrices of the different multigrid levels are cached during the solver setup phase to increase the efficiency of the overall solution process. The storage cost of these matrices is minimal and scales inversely in proportion to the number of MPI tasks. 
    \paragraph{Computational costs}\mbox{}\\
    The two main procedures that influence the overall computational cost of the proposed multigrid solution scheme are: \textit{i)} constructing and applying the additive Schwarz smoothers and \textit{ii)} solving the coarse problem. 
    \newline
    Construction of the AS smoothers following \eqref{eq:AdditiveSchwarz} entails the inversion of sub-matrices derived from basis function groups. The cost of this operation is dictated by the number and size of these sub-matrices. The number of sub-matrices increases linearly with the number of elements. The maximum size of a sub-matrix formed from the additive Schwarz groups, denoted by $m_{\text{max}}$, is proportional to the polynomial order $p$, the spatial dimension $d$, the manner in which the AS groups are selected (elementwise or patchwise selection) and the number of unknown field variables denoted by $n_f$. For a Poisson problem $n_f = 1$ while $n_f = 3$ for a three-dimensional linear elastic problem. Since the number of DOFs associated with the topological components for the tensor product and trunk space elements is known \cite{Duster2017}, it is possible to compute an upper bound for $m_{\text{max}}$ for uniform grids, see Table \ref{tab::upperbound} and Figure \ref{fig::mUniform}.
    \begin{table}[H]
        \small
\renewcommand{\arraystretch}{1.4}
    \centering
    \begin{tabular}{|c|c|}
    \toprule
        \textbf{tensor product space} & \textbf{trunk space} \\
    \midrule
        $n_f (n p + 1 )^3 $ &   $n_f (n+1)^2 ( 3 n p - 2n + 1 )$  \ for $p < 4$\\
                           &   $0.5 n_f (n+1)(3 n^2 p^2 - 9 n^2 p + 6 n p + 14 n^2 - 2 n + 2 )$ \ for $4 \leq p \leq 5$\\
                           &   $n_f (n^3 p^3 - 3 n^3 p^2 + 9 n^2 p^2 + 20 n^3 p - 9 n^2p + 18 (n p - n^3 + 2 n^2) + 6) $ \  for $p \geq 6$ \\
    \bottomrule
    \end{tabular}
        \caption{Maximum size of the additive Schwarz groups, $m_{\text{max}}$, in a uniform three-dimensional grid. $p$ represents the element polynomial order, $n_f$ the number of field variables in the problem and $n$ is a factor that is equal to one for elementwise blocks and equal to two when patchwise block selection is applied.}
        \label{tab::upperbound}
\end{table}
In the case of multi-level $hp$-grids, $m_{\text{max}}$ is not bounded and its size depends on the refinement level $k$ and the refinement pattern applied to the mesh elements. Figure \ref{fig::refUniform} shows the value of $m_{\text{max}}$ for the benchmark problem considered in Section \ref{sec::cubeSpheresEx}. From Table \ref{tab::upperbound} and Figure \ref{fig::boundsGraph} it is clear that the inversion of the patchwise additive Schwarz sub-matrices can become increasingly expensive for high polynomial orders and refinement levels. It is therefore important that optimized algorithms are used to perform these inversions. Our code framework makes use of distributed and shared memory parallelism to accelerate the construction of the additive Schwarz smoothers. For \enquote{small} sub-matrices, where $\mathbf{A}_i \in \mathbb{R}^{m \times m}$ and $m < 1000$, the built-in invert function of the \texttt{BOOST} library is used for the inversion. When $m$ exceeds 1000, the direct solver \texttt{Pardiso} \cite{Schenk2011} is used to invert $\mathbf{A}_i$ by solving an equation system with $m$ righthand side vectors. This operation can be written as $\mathbf{A}_i\mathbf{X} = \mathbf{B}$, where $\mathbf{B}$ is a $m \times m$ matrix whose columns are made up  of the unit vectors $\mathbf{e}_j$ with $j \in [0,m]$.  
The use of the \texttt{Epetra} package for the storage of the additive Schwarz smoothers $\mathbf{M}^{-1}_{\ell}$ allows our numerical code to make use of its optimized parallel multiplication kernels. This ensures that the smoothers can be applied in an efficient and scalable manner. We again leverage the \texttt{AztecOO} package in \texttt{Trilinos} for the solution of the coarse systems with $p=1$ and $k=0$. These systems are solved using the package's parallel CG solver and the additive Schwarz preconditioner published in \cite{Jomo2019}.  

\renewcommand{\picsDir}{examples/pics}
\renewcommand{\graphDir}{examples/graphs}

\section{Numerical examples}\label{sec::numerics}
The following section investigates the performance of the proposed hierarchical multigrid framework in a series of numerical examples. We use the multigrid V-cycle algorithm primarily as a preconditioner within the CG method but also show results in which the algorithm is utilized as a stand-alone solver. The quality of the solution obtained in the different examples is measured by monitoring the norm of the relative residual defined as
\begin{equation}
    \frac{\| \mathbf{r}_i \|_2}{\| \mathbf{b} \|_2} = \frac{\| \mathbf{b} - \mathbf{A} \mathbf{x}_i \|_2}{\| \mathbf{b} \|_2}\, 
\end{equation}
with the righthand side vector $\mathbf{b}$ and the terms $\mathbf{r}_i$ and $\mathbf{x_i}$ that denote the residual and approximation of the finest grid in the $i^{\textrm{th}}$ iteration, respectively. When the multigrid algorithm is used as a solver, the contraction  number $\rho_i$, defined as the quotient between two consecutive residual norms i.e.,
\begin{equation}\label{eq::contraction}
    \rho_i = \frac{\| \mathbf{r}_i \|_2}{ \| \mathbf{r}_{i -1} \|_2},
\end{equation}
is also used as a measure to judge the effectiveness of the solution scheme. $\rho_{max}$ is a measure that denotes the maximum contraction number (largest reduction factor) in a series of iterations. Only symmetric positive definite systems arising from problems in linear elasticity or the Poisson equation are considered in this section. A total of five pre- and post-smoothing steps are applied on every multigrid level $\ell \neq 0$ in each simulation. All computations are performed on the SuperMUC-NG system hosted at the Leibniz Supercomputing Center in Garching, Germany. The compute nodes used have an architecture comprising dual-socket Intel Xeon Scalable Platinum 8168 processors (Skylake). Each node has a total of 48 cores and 96GB of main memory. Furthermore, the following compiler and library versions are used: IntelMPI compiler version 19.0, Trilinos 12.12.1 and the GNU compiler 7.0. The most aggressive level of compiler optimization flags are chosen for the nodal architecture and include -O3 and -funroll-loops among others.

\subsection{Poisson problem on a square domain}\label{sec::squarePoisson}
\begin{figure}[H]
    \centering
    \subfloat[Geometry of the rotating square domain.]
    {
      \includegraphics[width=0.40\textwidth]{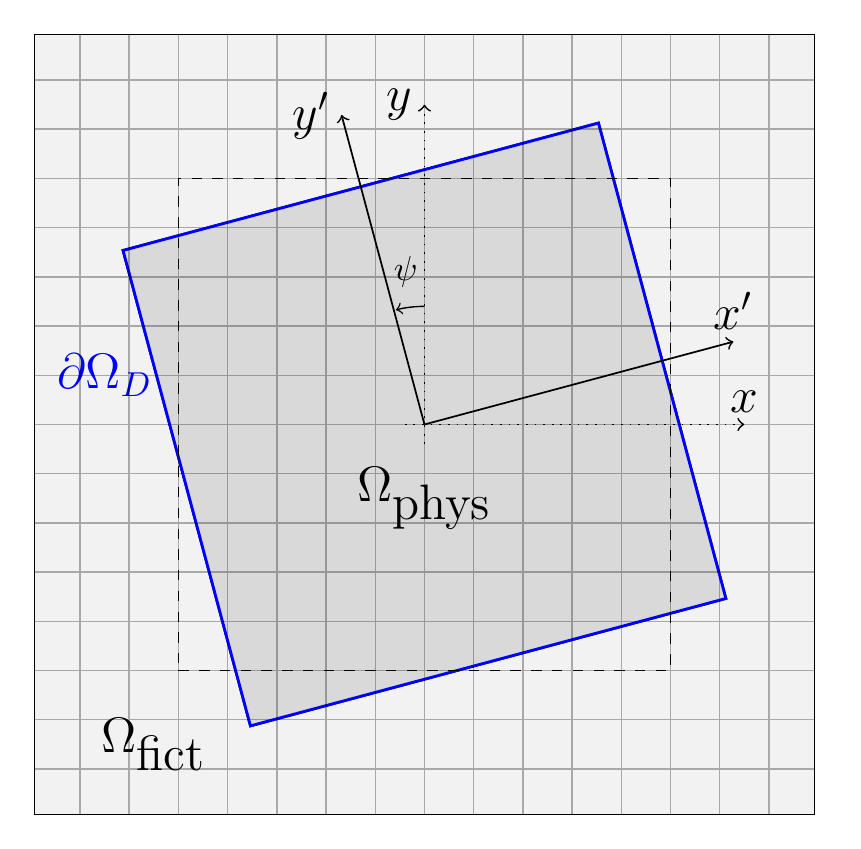}
    \label{fig::PoissonSetup}
    }
    \hspace{5mm}
    \subfloat[Prescribed solution.]
    {
      \includegraphics[width=0.40\textwidth]{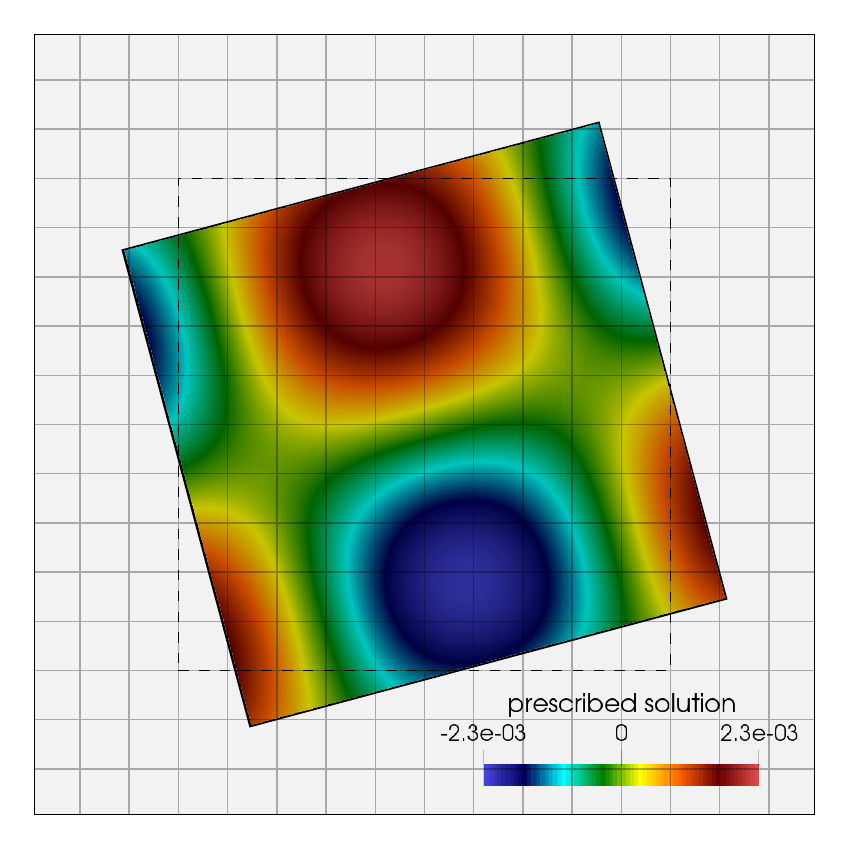}
    }
    \caption{Rotating Poisson problem with a manufactured solution.}
    \label{fig::PoissonExample}
\end{figure}

The first numerical example investigates the performance of the proposed $hp$-multigrid algorithm as a stand-alone solver and as a preconditioner within a CG scheme. It aims to demonstrate the effectiveness of the additive Schwarz technique when used as a smoother for finite cell problems on uniform meshes. We consider a Poisson problem posed on a square domain of unit length in which the physical domain is rotated about the origin with respect to a fixed background grid as shown in Figure \ref{fig::PoissonExample}. This setup results in different cut scenarios when the angle of rotation $\psi$ is varied. The method of manufactured solutions is applied in this numerical test and a temperature field $\bar{u}$ is chosen such that

\begin{equation}  
  \bar{u} = \frac{1}{2 \kappa a^2} \cos(a\,x^\prime)\sin(a\,y^\prime), 
\end{equation} 
 with $a = \frac{3}{2}\pi$ and $x^\prime,y^\prime$ denoting the coordinates of the rotated coordinate system defined at the center of the domain. A source term $s$ is derived from the manufactured solution following the Poisson equation $-\kappa \Delta u = s$ with $ s = \cos(a\,x^\prime)\sin(a\,y^\prime)$ and applied in $\Omega_{\text{phys}}$. Dirichlet boundary conditions are prescribed on all edges of the square domain such that $u = \bar{u}$ on $\Gamma_D$. These constraints are enforced using the penalty method with $\beta = 10^4$. A value of $\alpha = 10^{-8}$ is applied in $\Omega_{\textrm{fict}}$ and a value of $\kappa = 10$ is used in $\Omega_{\text{phys}}$. 


\subsubsection{Convergence behavior in the mesh-fitting case}\label{sec::poissonmeshfitting}
We first analyze the performance of the hierarchical multigrid scheme in the mesh-fitting case where $\psi = 0$. We consider uniform grids using high-order elements with different polynomial orders with $p \in [ 2, \dots, 5]$ and mesh sizes $h  = \lbrace{ \frac{1}{8}, \frac{1}{16}, \frac{1}{32}, \frac{1}{64} \rbrace}$. Moreover, the effectiveness of four different smoothers is investigated: \textit{i)} Jacobi smoothing, \textit{ii)} Gauss-Seidel smoothing, \textit{iii)} additive Schwarz smoothing with elementwise blocks and \textit{iv)} additive Schwarz smoothing with patchwise blocks. In the study at hand, the convergence behavior of the hierarchical multigrid approach as a stand-alone solver and a preconditioner is analyzed for the different combinations of smoother, polynomial orders and mesh sizes. Each solver is terminated when the value of the relative residual is below
$10^{-9}$ or when the number of iterations exceeds 500. The $*$ symbol is used to represent scenarios in which the solver did not converge within 500 iterations. The results of this study are summarized in Table \ref{tab:testcaseA}. 

Table \ref{tab:testcaseA} shows the convergence behavior of the multigrid algorithm as a solver and preconditioner for $\psi = 0^{\circ}$. In this boundary-conforming example, all smoothers achieve convergence rates that are independent of the mesh size $h$. These results indicate the correct implementation of the multigrid scheme. These convergence rates are, however, dependent on the polynomial order for the Jacobi, Gauss-Seidel and elementwise additive Schwarz smoothers as shown in Table \ref{tab:testcaseA} and indicated by the maximum contraction of the residual for a mesh with $h=\frac{1}{32}$ recorded in Table \ref{tab::squareCont}. These smoothers appear to perform better for odd polynomial orders than for even orders. A similar odd-even pattern is reported in \citep{Babuska1989} and may be due to the asymmetry of the manufactured solution. Convergence rates independent of the polynomial order $p$ are obtained using the patchwise additive Schwarz smoother.

\subsubsection{Convergence behavior in immersed configurations}
The multigrid algorithm is now used in an immersed setting, in which the domain $\Omega_{\textrm{phys}}$ is rotated about the origin as shown in Figure \ref{fig::PoissonSetup}. In analogy to Section \ref{sec::squarePoisson}, we compare the performance of different smoothers for $\psi=30^{\circ}$ and summarize the results in Table \ref{tab:testcaseB}. The effect of the ill-conditioning due to the cut cells is clearly seen in this example as the standard smoothers, i.e.\ the Jacobi and Gauss-Seidel methods, fail to improve the conditioning of the linear systems and are characterized by poor convergence of the multigrid solution techniques applied in this study. The additive Schwarz smoothers do a better job of detecting almost linear dependent functions as shown in Table \ref{tab:testcaseB}. Note that the elementwise smoother may contain a few small modes that cause slow convergence when multigrid is utilized as a stand-alone solver. These modes, however, only result in a few additional CG iterations, when the multigrid algorithm is used as a preconditioner. The patchwise smoother robustly deals with all small modes due to cut cells and results in a convergence behavior that is independent of $h$ and $p$. Table \ref{tab::squareImmersedCont} presents the maximum contraction number for different polynomial orders of a mesh with $h = \frac{1}{32}$.

\begin{table}[H]
    \small
    \paragraph{Rotating Poisson problem: Convergence study for $\psi = 0^{\circ}$}
\centering
\subfloat
{
\begin{tabularx}{0.5\textwidth}{|c| *{4}{Y|}}
\toprule
   & \multicolumn{4}{c|}{Multigrid as a solver}\\
	\cline{2-5}  & \multicolumn{4}{c|}{Jacobi smoother (5,5)}\\
\toprule
    \diagbox{$p$}{$h$} & $\frac{1}{8}$ &  $\frac{1}{16}$ & $\frac{1}{32}$ & $\frac{1}{64}$ \\
\toprule
	2 &   23  &   25  &   25   &   24   \\
	3 &   17  &   16  &   15   &   15   \\
	4 &   31  &   33  &   35   &   37   \\
	5 &   24  &   24  &   26   &   27   \\
\bottomrule
\end{tabularx}
}
\subfloat
{
\begin{tabularx}{0.5\textwidth}{|c| *{4}{Y|}}
\toprule
   & \multicolumn{4}{c|}{CG with multigrid preconditioner}\\
	\cline{2-5}  & \multicolumn{4}{c|}{Jacobi smoother (5,5)}\\
\toprule
    \diagbox{$p$}{$h$} & $\frac{1}{8}$ &  $\frac{1}{16}$ & $\frac{1}{32}$ & $\frac{1}{64}$ \\
\toprule
	2 &  10   &  10   &   10   &    10  \\
	3 &   8   &   8   &    8   &     8  \\
	4 &  13   &  13   &   13   &    13  \\
	5 &  10   &  11   &   11   &    11  \\
\bottomrule
\end{tabularx}
}
\hfill
\subfloat
{
\begin{tabularx}{0.5\textwidth}{|c| *{4}{Y|}}
\toprule
   & \multicolumn{4}{c|}{Multigrid as a solver}\\
	\cline{2-5}  & \multicolumn{4}{c|}{Gauss-Seidel smoother (5,5)}\\
\toprule
    \diagbox{$p$}{$h$} & $\frac{1}{8}$ &  $\frac{1}{16}$ & $\frac{1}{32}$ & $\frac{1}{64}$ \\
\toprule
	2 &   7   &    7  &    7   &    7    \\
	3 &   6   &    6  &    6   &    6    \\
	4 &  13   &   13  &   14   &   14     \\
	5 &  10   &   10  &   10   &   10     \\
\bottomrule
\end{tabularx}
}
\subfloat
{
\begin{tabularx}{0.5\textwidth}{|c| *{4}{Y|}}
\toprule
   & \multicolumn{4}{c|}{CG with multigrid preconditioner}\\
	\cline{2-5}  & \multicolumn{4}{c|}{symm. Gauss-Seidel smoother (5,5)}\\
\toprule
    \diagbox{$p$}{$h$} & $\frac{1}{8}$ &  $\frac{1}{16}$ & $\frac{1}{32}$ & $\frac{1}{64}$ \\
\toprule
	2 &   6   &    6  &    6   &    6   \\
	3 &   5   &    5  &    5   &    5   \\
	4 &   8   &    8  &    8   &    8   \\
	5 &   8   &    7  &    7   &    7   \\
\bottomrule
\end{tabularx}
}
\hfill
\subfloat
{
\begin{tabularx}{0.5\textwidth}{|c| *{4}{Y|}}
\toprule
   & \multicolumn{4}{c|}{Multigrid as a solver}\\
	\cline{2-5}  & \multicolumn{4}{c|}{elementwise AS smoother (5,5)}\\
\toprule
    \diagbox{$p$}{$h$} & $\frac{1}{8}$ &  $\frac{1}{16}$ & $\frac{1}{32}$ & $\frac{1}{64}$ \\
\toprule
	2 &   12    &   13    &   13     &   13    \\
	3 &   10    &   10    &   10     &    9    \\
	4 &   17    &   17    &   17     &   17    \\
	5 &   14    &   12    &   13     &   13    \\
\bottomrule
\end{tabularx}
}
\subfloat
{
\begin{tabularx}{0.5\textwidth}{|c| *{4}{Y|}}
\toprule
   & \multicolumn{4}{c|}{CG with multigrid preconditioner}\\
	\cline{2-5}  & \multicolumn{4}{c|}{elementwise AS smoother (5,5)}\\
\toprule
    \diagbox{$p$}{$h$} & $\frac{1}{8}$ &  $\frac{1}{16}$ & $\frac{1}{32}$ & $\frac{1}{64}$ \\
\toprule
	2 &  7   &   8  &   8  &   8   \\
	3 &  7   &   7  &   7  &   7   \\
	4 &  9   &   9  &   9  &   9     \\
	5 &  8   &   8  &   8  &   8     \\
\bottomrule
\end{tabularx}
}
\hfill
\subfloat
{
\begin{tabularx}{0.5\textwidth}{|c| *{4}{Y|}}
\toprule
   & \multicolumn{4}{c|}{Multigrid as a solver}\\
	\cline{2-5}  & \multicolumn{4}{c|}{patchwise AS smoother (5,5)}\\
\toprule
    \diagbox{$p$}{$h$} & $\frac{1}{8}$ &  $\frac{1}{16}$ & $\frac{1}{32}$ & $\frac{1}{64}$ \\
\toprule
	2 &   8   &   8   &    8   &    7   \\
	3 &   7   &   6   &    6   &    6   \\
	4 &   7   &   7   &    7   &    6   \\
	5 &   6   &   6   &    5   &    5   \\
\bottomrule
\end{tabularx}
}
\subfloat
{
\begin{tabularx}{0.5\textwidth}{|c| *{4}{Y|}}
\toprule
   & \multicolumn{4}{c|}{CG with multigrid preconditioner}\\
	\cline{2-5}  & \multicolumn{4}{c|}{patchwise AS smoother (5,5)}\\
\toprule
    \diagbox{$p$}{$h$} & $\frac{1}{8}$ &  $\frac{1}{16}$ & $\frac{1}{32}$ & $\frac{1}{64}$ \\
\toprule
	2 &  6  &   6  &   6  &   6   \\
	3 &  6  &   5  &   5  &   5   \\
	4 &  6  &   6  &   6  &   6   \\
	5 &  5  &   5  &   5  &   5   \\
\bottomrule
\end{tabularx}
}
\caption{Convergence study for $\psi = 0$ comparing the performance of different smoothers for varying element sizes and polynomial orders. The figures in the tables represent the number of iterations required to reach a tolerance in the relative residual of $10^{-9}$. Five pre- and post-smoothing steps are performed per V-cycle on each multigrid level with $\ell \neq 0$.}
\label{tab:testcaseA}
\end{table}

\begin{table}[H]
    \centering
      \begin{tabularx}{0.6\textwidth}{|c| *{4}{Y|}}
      \toprule
          \diagbox{Smoother}{$p$} & 2 &  3 & 4 & 5 \\
      \toprule
      	Jacobi         &  .478  &  .391  &   .683  &   .617  \\
      	Gauss-Seidel   &  .125  &  .095  &   .362  &   .255  \\
      	Elementwise AS &  .309  &  .252  &   .430  &   .414  \\
      	Patchwise AS   &  .162  &  .164  &   .140  &   .162  \\
      \bottomrule
          \end{tabularx}
          \caption{Maximum contraction number $\rho_{\text{max}}$ of the $p$-multigrid solver for $h = \frac{1}{32}$ and $\psi = 0^{\circ}$.}
          \label{tab::squareCont}
\end{table}

\begin{table}[H]
    \paragraph{Rotating Poisson problem: Convergence study for $\psi = 30^{\circ}$}
    \small
\centering
\subfloat
{
\begin{tabularx}{0.5\textwidth}{|c| *{4}{Y|}}
\toprule
   & \multicolumn{4}{c|}{Multigrid as a solver}\\
	\cline{2-5}  & \multicolumn{4}{c|}{Jacobi smoother (5,5)}\\
\toprule
    \diagbox{$p$}{$h$} & $\frac{1}{8}$ &  $\frac{1}{16}$ & $\frac{1}{32}$ & $\frac{1}{64}$ \\
\toprule
 	2 &    *  &    *  &    *   &    *   \\
 	3 &    *  &    *  &    *   &    *   \\
 	4 &    *  &    *  &    *   &    *   \\
 	5 &    *  &    *  &    *   &    *   \\
\bottomrule
\end{tabularx}
}
\subfloat
{
\begin{tabularx}{0.5\textwidth}{|c| *{4}{Y|}}
\toprule
   & \multicolumn{4}{c|}{CG with multigrid preconditioner}\\
	\cline{2-5}  & \multicolumn{4}{c|}{Jacobi smoother (5,5)}\\
\toprule
    \diagbox{$p$}{$h$} & $\frac{1}{8}$ &  $\frac{1}{16}$ & $\frac{1}{32}$ & $\frac{1}{64}$ \\
\toprule
    2 &  133   &  162   &   201   &   180  \\
	3 &   *   &   *   &    *   &     *  \\
	4 &   *   &   *   &    *   &     *  \\
	5 &   *   &   *   &    *   &     *  \\
\bottomrule
\end{tabularx}
}
\hfill
\subfloat
{
\begin{tabularx}{0.5\textwidth}{|c| *{4}{Y|}}
\toprule
   & \multicolumn{4}{c|}{Multigrid as a solver}\\
	\cline{2-5}  & \multicolumn{4}{c|}{Gauss-Seidel smoother (5,5)}\\
\toprule
    \diagbox{$p$}{$h$} & $\frac{1}{8}$ &  $\frac{1}{16}$ & $\frac{1}{32}$ & $\frac{1}{64}$ \\
\toprule
    2 &    *  &   *   &   *    &   *    \\
	3 &    *  &   *   &   *    &   *    \\
	4 &    *  &   *   &   *    &   *    \\
	5 &    *  &   *   &   *    &   *    \\
\bottomrule
\end{tabularx}
}
\subfloat
{
\begin{tabularx}{0.5\textwidth}{|c| *{4}{Y|}}
\toprule
   & \multicolumn{4}{c|}{CG with multigrid preconditioner}\\
	\cline{2-5}  & \multicolumn{4}{c|}{symm. Gauss-Seidel smoother (5,5)}\\
\toprule
    \diagbox{$p$}{$h$} & $\frac{1}{8}$ &  $\frac{1}{16}$ & $\frac{1}{32}$ & $\frac{1}{64}$ \\
\toprule
    2 &  8   &  139  &  240   &    279 \\
	3 &   *   &   *   &    *   &     *  \\
	4 &   *   &   *   &    *   &     *  \\
	5 &   *   &   *   &    *   &     *  \\
\bottomrule
\end{tabularx}
}
\hfill
\subfloat
{
\begin{tabularx}{0.5\textwidth}{|c| *{4}{Y|}}
\toprule
   & \multicolumn{4}{c|}{Multigrid as a solver}\\
	\cline{2-5}  & \multicolumn{4}{c|}{elementwise AS smoother (5,5)}\\
\toprule
    \diagbox{$p$}{$h$} & $\frac{1}{8}$ &  $\frac{1}{16}$ & $\frac{1}{32}$ & $\frac{1}{64}$ \\
\toprule
	2 &   19  &   15  &   13   &   12   \\
	3 &   19  &   14  &   12   &   11   \\
	4 &   20  &   18  &   17   &   17   \\
	5 &   20  &   16  &   15   &   15   \\
\bottomrule
\end{tabularx}
}
\subfloat
{
\begin{tabularx}{0.5\textwidth}{|c| *{4}{Y|}}
\toprule
   & \multicolumn{4}{c|}{CG with multigrid preconditioner}\\
	\cline{2-5}  & \multicolumn{4}{c|}{elementwise AS smoother (5,5)}\\
\toprule
    \diagbox{$p$}{$h$} & $\frac{1}{8}$ &  $\frac{1}{16}$ & $\frac{1}{32}$ & $\frac{1}{64}$ \\
\toprule
	2 &   9   &   9   &    9   &     8  \\
	3 &  11   &  10   &    9   &     8  \\
	4 &  11   &  11   &   10   &    10  \\
	5 &  12   &  10   &    9   &     9  \\
\bottomrule
\end{tabularx}
}
\hfill
\subfloat
{
\begin{tabularx}{0.5\textwidth}{|c| *{4}{Y|}}
\toprule
   & \multicolumn{4}{c|}{Multigrid as a solver}\\
	\cline{2-5}  & \multicolumn{4}{c|}{patchwise AS smoother (5,5)}\\
\toprule
    \diagbox{$p$}{$h$} & $\frac{1}{8}$ &  $\frac{1}{16}$ & $\frac{1}{32}$ & $\frac{1}{64}$ \\
\toprule
	2 &    8  &    8  &    7   &    6   \\
	3 &    6  &    7  &    6   &    5   \\
	4 &    5  &    6  &    5   &    5   \\
	5 &    5  &    5  &    5   &    4   \\
\bottomrule
\end{tabularx}
}
\subfloat
{
\begin{tabularx}{0.5\textwidth}{|c| *{4}{Y|}}
\toprule
   & \multicolumn{4}{c|}{CG with multigrid preconditioner}\\
	\cline{2-5}  & \multicolumn{4}{c|}{patchwise AS smoother (5,5)}\\
\toprule
    \diagbox{$p$}{$h$} & $\frac{1}{8}$ &  $\frac{1}{16}$ & $\frac{1}{32}$ & $\frac{1}{64}$ \\
\toprule
	2 &   7   &   6   &    6   &     6  \\
	3 &   6   &   6   &    6   &     5  \\
	4 &   5   &   6   &    5   &     5  \\
	5 &   5   &   5   &    5   &     5  \\
\bottomrule
\end{tabularx}
}
\caption{Convergence study for $\psi = 30^{\circ}$ comparing the performance of different smoothers for varying element sizes and polynomial orders. The figures in the tables represent the number of iterations required by the solver while the symbol * indicates that the solver did not converge to a tolerance of $10^{-9}$ within 500 iterations. Five pre- and post-smoothing steps are performed per V-cycle on each multigrid level with $\ell \neq 0$}
\label{tab:testcaseB}
\end{table}

\begin{table}[H]
    \centering
      \begin{tabularx}{0.6\textwidth}{|c| *{4}{Y|}}
      \toprule
          \diagbox{Smoother}{$p$} & 2 &  3 & 4 & 5 \\
      \toprule
      	Jacobi         &  .998  &  .998  &   .998  &   .998  \\
      	Gauss-Seidel   &  .996  &  div.   &   div.   &   div.   \\
      	Elementwise AS &  .356  &  .355  &   .521  &   .484  \\
      	Patchwise AS   &  .142  &  .117  &   .100  &   .108  \\
      \bottomrule
          \end{tabularx}
          \caption{Maximum contraction number $\rho_{\text{max}}$ of the $p$-multigrid solver for $h = \frac{1}{32}$ and $\psi = 30^{\circ}$. The symbol div. indicates that the solver diverged. }
          \label{tab::squareImmersedCont}
\end{table}
\subsection{Perforated linear elastic plate}
\begin{figure}[H]
 \centering
 \subfloat[Perforated linear elastic plate.]
 {
  \includegraphics[width=0.32\textwidth]{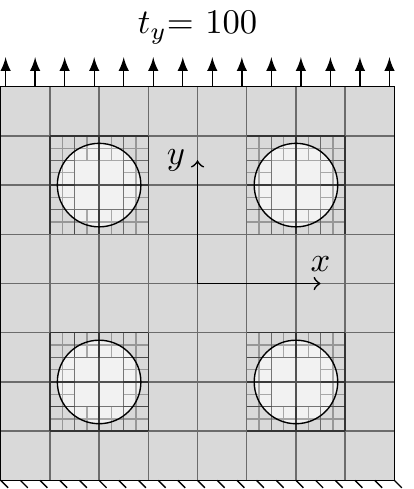}
  \label{fig::perforatedPlateGeom}
 }
 \hspace{2cm}
 \subfloat[Von Mises stress.]
 {
  \includegraphics[width=0.32\textwidth]{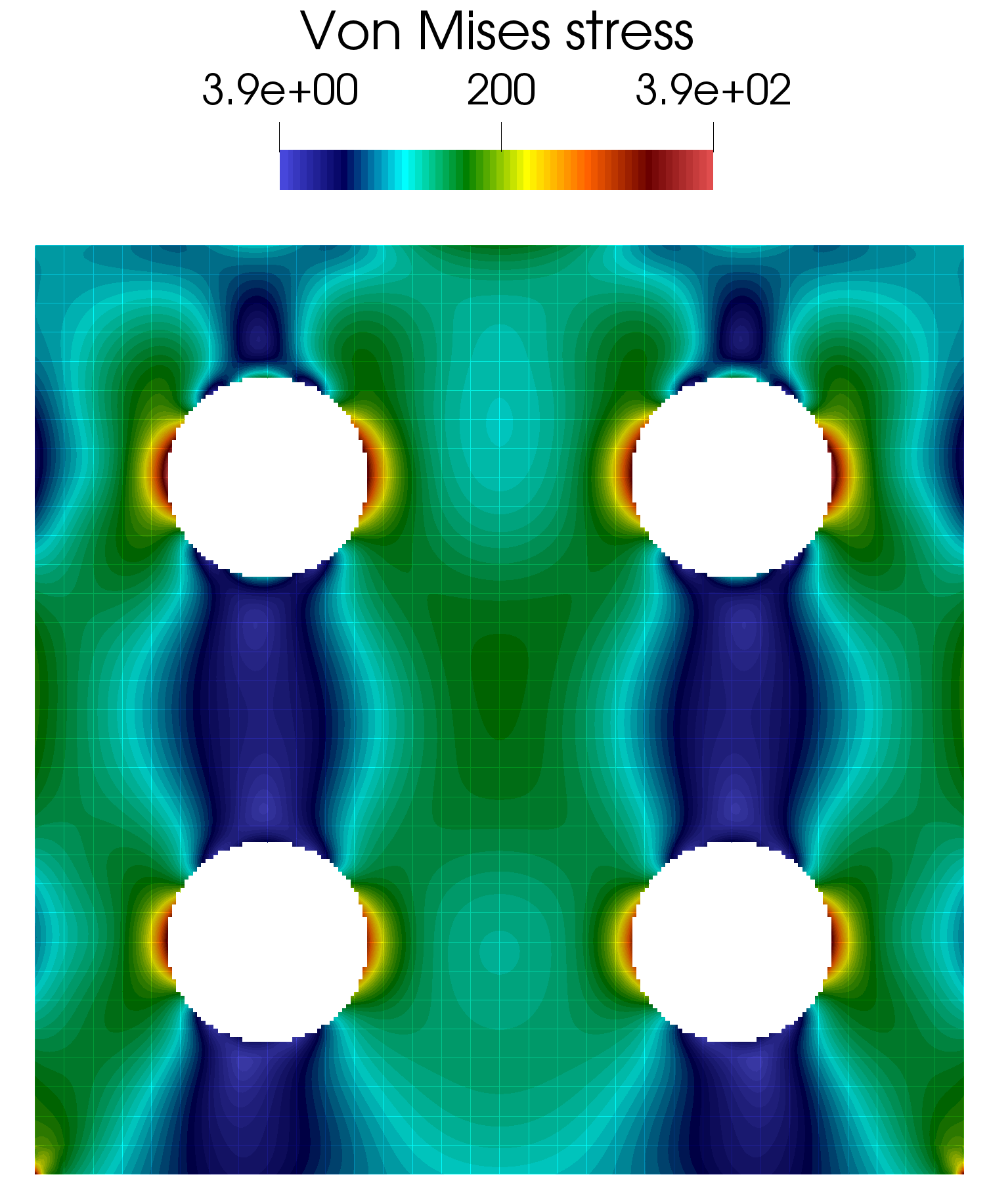}
 }
 \caption{Setup and von Mises stress of the linear elastic perforated plate example.}
\end{figure}

The previous numerical example established the suitability of using a multigrid preconditioner in conjunction with additive Schwarz smoothers for FCM problems arising from the Poisson equation. The performance of this preconditioner is now investigated in the context of linear elasticity. To this end, a square domain with a length of $l$ = 4 is subjected to a tensional force on one end and fully clamped on the other end as depicted in Figure \ref{fig::perforatedPlateGeom}. The square domain has four circular cavities with a radius of $0.3\sqrt{2}$ and is characterized by an elastic modulus $E = 2.069 \cdot 10^{5}$ MPa and a Poisson's ratio $\nu$ = 0.29. The finite cell method is applied in this example with $\alpha= 10^{-8}$ in $\Omega_{\text{fict}}$ and a penalty parameter $\beta=10^{8}$. The number of elements per direction is chosen such that $h \in \lbrace{ \frac{1}{8}, \frac{1}{16}, \frac{1}{32}, \frac{1}{64} \rbrace}$.

Using the setup shown in Figure \ref{fig::perforatedPlateGeom}, numerical studies are carried out that investigate the convergence of a CG solver when the presented additive Schwarz techniques are employed as preconditioners or as smoothers within an $hp$-multigrid preconditioner. The study is conducted on multi-level $hp$-refined grids with a fixed polynomial order and varying levels of refinement. Elements that are intersected by the immersed boundary are refined recursively to a predefined depth as shown in Figure \ref{fig::perforatedPlateGeom}.

Table \ref{tab::perPlateRef} shows the number of iterations required by the different preconditioner and smoother configurations for multilevel $hp$-grids with quadratic elements and varying levels of refinement. When the elementwise additive Schwarz blocks are utilized as both smoothers and preconditioners, increasing the refinement depth leads to an increase in the number of iterations. The multigrid algorithm that uses these elementwise blocks for smoothing does not achieve convergence rates that are independent of the mesh parameter $h$. This behavior is attributed to the fact that certain small modes can remain untreated by the elementwise blocks leading to slow convergence. In contrast, the patchwise blocks yield convergence rates that are independent of the refinement depth employed. Furthermore, the multigrid algorithm that employs a patchwise smoothing approach achieves mesh-independent convergence rates.  

\begin{table}[H]
\centering
    \subfloat
{
\begin{tabularx}{0.5\textwidth}{|c| *{4}{Y|}}
\toprule
   & \multicolumn{4}{c|}{CG with elementwise AS}\\
	\cline{2-5}  & \multicolumn{4}{c|}{preconditioner}\\
\toprule
    \diagbox{$k$}{$h$} & $\frac{1}{8}$ & $\frac{1}{16}$ & $\frac{1}{32}$ & $\frac{1}{64}$ \\
\toprule
    0 & 60 & 76 &   135 &  262      \\
	1 & 81  &  220 &  234 & 1192      \\
	2 & 183 &  367 &  807 & 2114      \\
	3 & 294 &  483 & 1803 & 2566      \\
\bottomrule
\end{tabularx}
}
    \subfloat
{
\begin{tabularx}{0.5\textwidth}{|c| *{4}{Y|}}
\toprule
   & \multicolumn{4}{c|}{CG with multigrid preconditioner}\\
	\cline{2-5}  & \multicolumn{4}{c|}{elementwise AS smoother (5,5)}\\
\toprule
    \diagbox{$k$}{$h$} & $\frac{1}{8}$ & $\frac{1}{16}$ & $\frac{1}{32}$ & $\frac{1}{64}$ \\
\toprule
	0 & 12 &    9  &   9  &   9   \\
	1 & 15 &   30  &  29  &  79   \\
	2 & 31 &   46  &  63  & 115   \\
	3 & 50 &   67  & 102  & 114   \\
\bottomrule
\end{tabularx}
}
\hfill
    \subfloat
{
\begin{tabularx}{0.5\textwidth}{|c| *{4}{Y|}}
\toprule
   & \multicolumn{4}{c|}{CG with patchwise AS }\\
	\cline{2-5}  & \multicolumn{4}{c|}{preconditioner}\\
\toprule
    \diagbox{$k$}{$h$} & $\frac{1}{8}$ & $\frac{1}{16}$ & $\frac{1}{32}$ & $\frac{1}{64}$ \\
\toprule
	0 &  43  & 64 & 109 & 210 \\
	1 &  43  & 64 & 109 & 210 \\
	2 &  43  & 65 & 109 & 210 \\
	3 &  43  & 64 & 109 & 210 \\
\bottomrule
\end{tabularx}
}
    \subfloat
{
\begin{tabularx}{0.5\textwidth}{|c| *{4}{Y|}}
\toprule
   & \multicolumn{4}{c|}{CG with multigrid preconditioner}\\
	\cline{2-5}  & \multicolumn{4}{c|}{patchwise AS smoother (5,5)}\\
\toprule
    \diagbox{$k$}{$h$} & $\frac{1}{8}$ & $\frac{1}{16}$ & $\frac{1}{32}$ & $\frac{1}{64}$ \\
\toprule
	2 &  7 &  7 &  7 &  7     \\
	3 &  7 &  7 &  7 &  7     \\
	4 &  7 &  7 &  6 &  6      \\
	5 &  6 &  6 &  6 &  6      \\
\bottomrule
\end{tabularx}
}
    \caption{Perforated plate example: Convergence behavior of a CG solver for four levels of refinement and four different preconditioners. The study considers multi-level $hp$-refined finite cell meshes with varying resolutions and a fixed polynomial order of $p=2$. }
    \label{tab::perPlateRef}
\end{table}

\subsection{Cube with spherical cavities}\label{sec::cubeSpheresEx}
\begin{figure}[H]
 \begin{center}
   \subfloat[Cube geometry.]
     {
       \label{fig::cubeWithCavities}
       \includegraphics[width=0.32\textwidth]{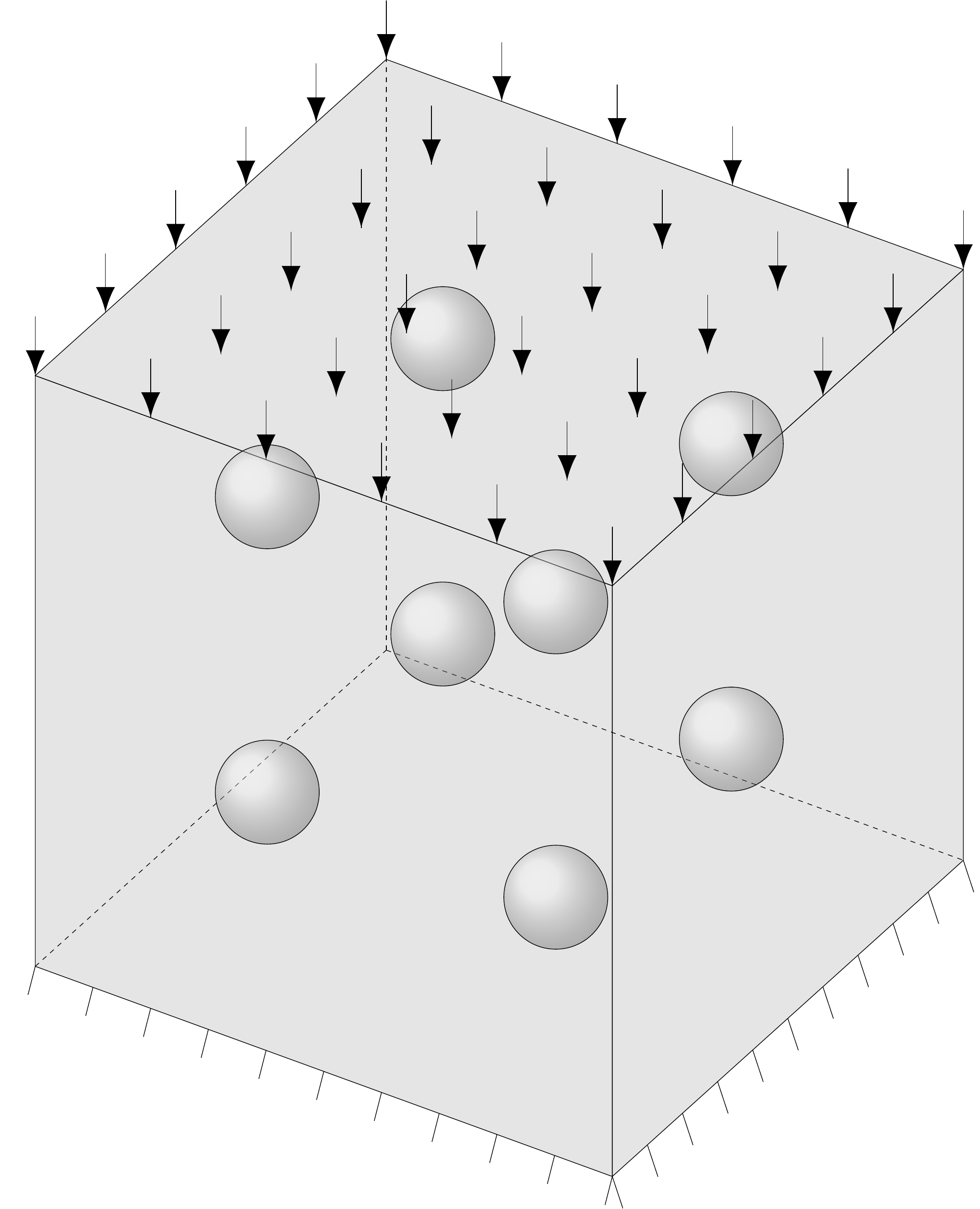}%
     }%
     \hspace{4mm}
   \subfloat[Von Mises stress.]
     {
       \includegraphics[width=0.48\textwidth]{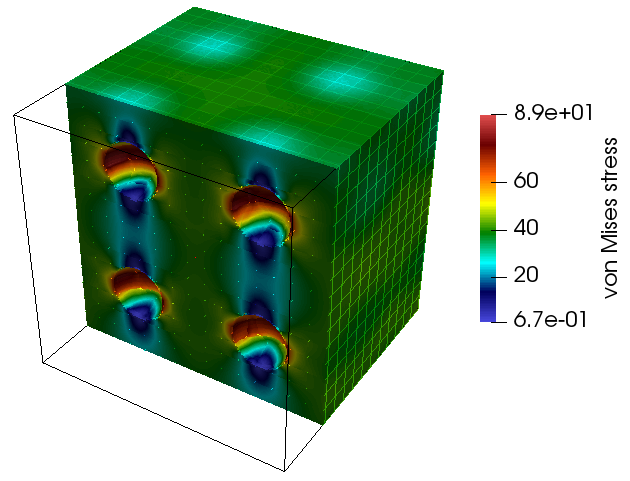}%
     }%
 \caption{Cube with spherical cavities example.}
 \end{center}
\end{figure}

To assess the performance of the proposed $hp$-multigrid approach in a three-dimensional setting, we now consider a simple example consisting of a linear elastic cube of unit length with spherical cavities subject to compressional loading. The cube has the same material properties as the perforated plate in the previous example and the values of $\alpha = 10^{-8}$ and $\beta = 10^{10}$ are chosen. The radii of the cavities are again $r = 0.3\sqrt{2}$. A homogeneous pressure load $P$ = 100 N/$\text{mm}^2$ is applied on the upper surface of the cube as shown in Figure \ref{fig::cubeWithCavities}. 

\subsection*{Influence of the polynomial order}
\begin{figure}[H]
 \begin{center}
   \subfloat[Relative residual for $p=2$.]
    {
       \includegraphics[width=0.48\textwidth]{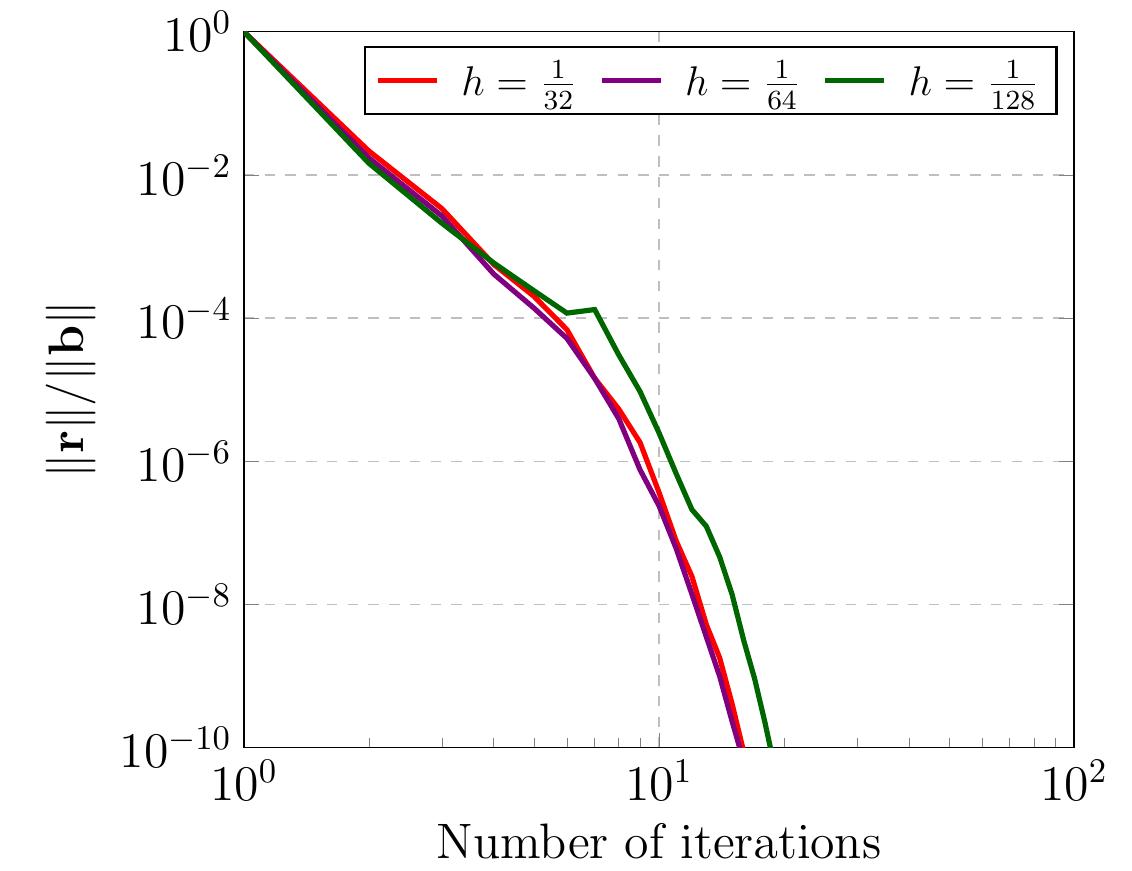}%
     }%
   \subfloat[Relative residual for $p=3$.]
     {
       \includegraphics[width=0.48\textwidth]{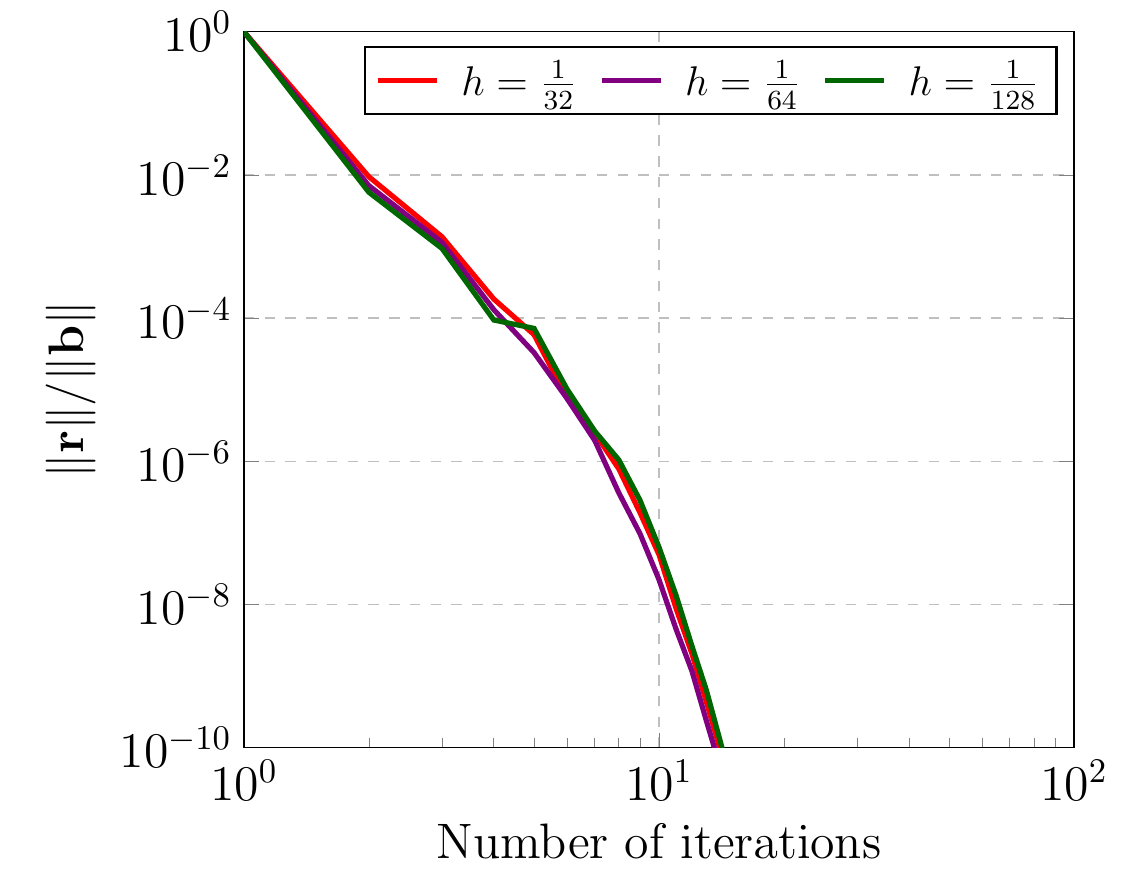}%
     }%
     \hfill
   \subfloat[Relative residual for $p=4$.]
     {
       \includegraphics[width=0.48\textwidth]{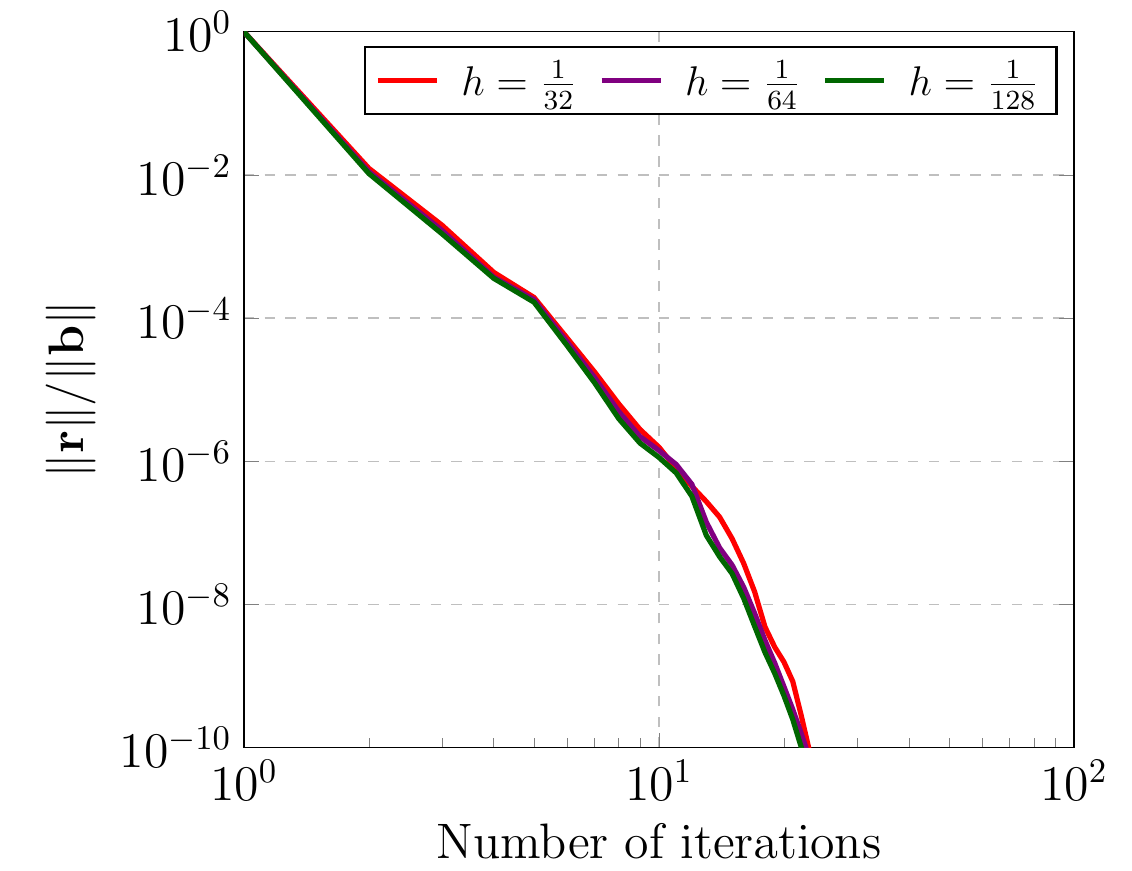}%
     }%
   \subfloat[Relative residual for $p=5$.]
     {
       \includegraphics[width=0.48\textwidth]{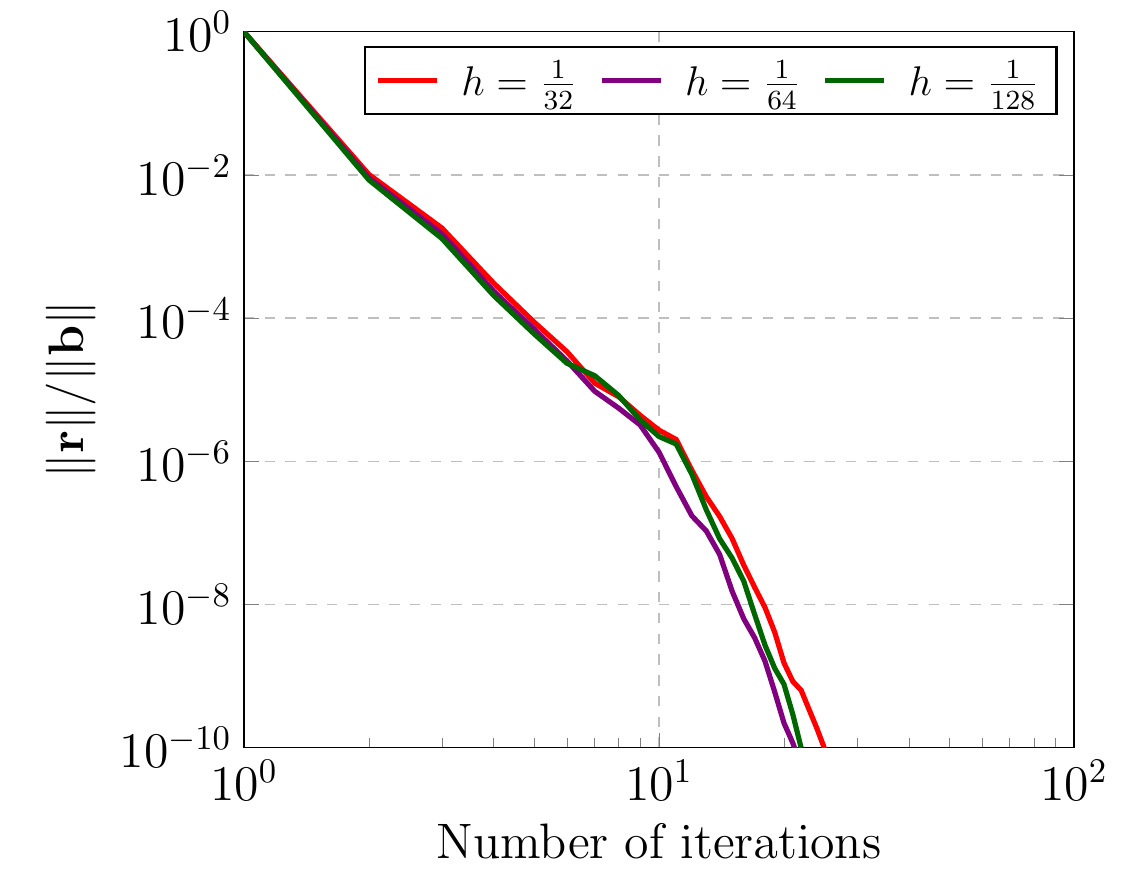}%
     }%
   \hfill
 \caption{Convergence of the CG solver with a $p$-multigrid preconditioner and elementwise additive Schwarz smoothing for hexahedral trunk space elements.}
     \label{fig::Cube4HolesConv}
 \end{center}
\end{figure}

The effect of the element polynomial order on the convergence of a CG solver with a $p$-multigrid preconditioner utilizing elementwise additive Schwarz smoothing is the subject of the following study. To this end, a sequence of meshes with varying element sizes is analyzed. The mesh size is chosen such that $h = \lbrace{ \frac{1}{32}, \frac{1}{64}, \frac{1}{128} \rbrace}$. The number of unknowns ranges from $417\,339$ for the coarsest $p=2$ mesh to $135\,864\,459$ for the finest mesh with $p=5$. Similar to the two-dimensional benchmark cases, the proposed $p$-multigrid approach leads to convergence rates that are independent of the mesh size $h$ as shown in Figure \ref{fig::Cube4HolesConv}. The number of iterations ranges between 16 and 19 for $p=2$, between 14 and 15  for $p=3$, between 22 and 23 for $p=4$ and between 22 and 26 for $p=5$.  



\begin{figure}[H]
 \centering
    \subfloat[CPU timings for $h=\frac{1}{32}$. All simulations are run on 96 cores.]
    {
     \label{fig::Cube4HolesCosth32}
       \includegraphics[width=0.48\textwidth]{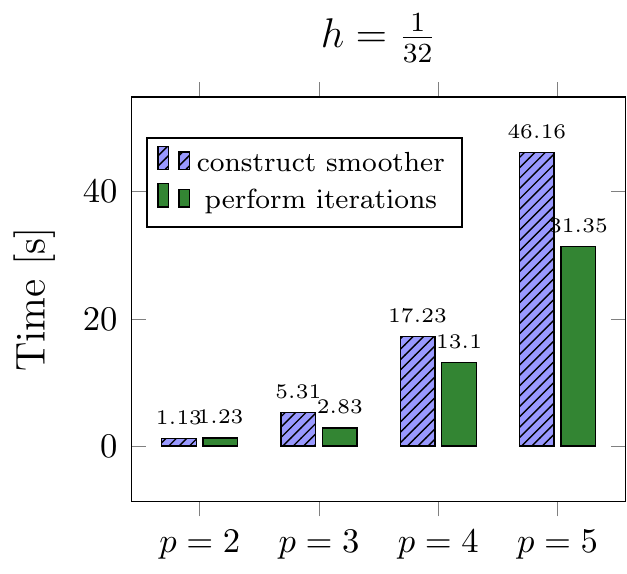}%
    }
    \hfill
    \subfloat[CPU timings for $h=\frac{1}{64}$. All simulations are run on 768 cores.]
    {
     \label{fig::Cube4HolesCosth64}
       \includegraphics[width=0.49\textwidth]{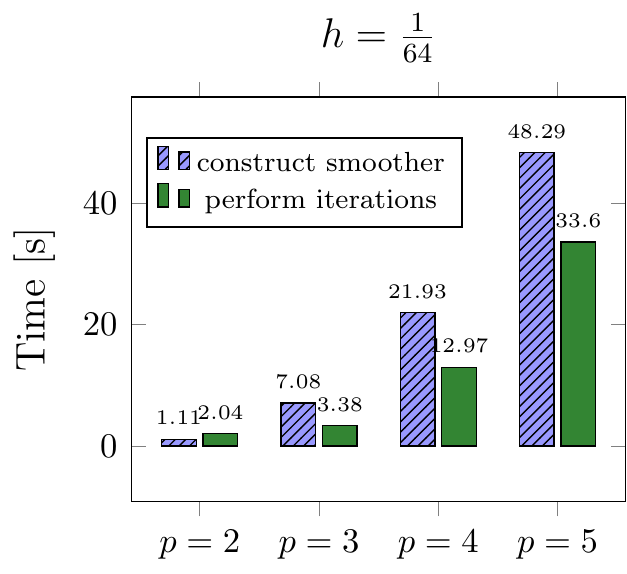}%
    }
 \caption{Execution time for the CG solver with a $p$-multigrid preconditioner.}
     \label{fig::Cube4HolesUniformTime}
\end{figure}

Figure \ref{fig::Cube4HolesUniformTime} shows the computational cost of the CG solver with a $p$-multigrid focusing on the time spent to construct the smoother and to perform the CG iterations. The latter procedure includes the application of the smoother and the solution of the coarse problem. From Figures \ref{fig::Cube4HolesCosth32} and \ref{fig::Cube4HolesCosth64} one can see that the computational cost of the solver increases for higher polynomial orders.

\subsection*{Influence of the refinement level}
\begin{figure}[H]
 \centering
    \subfloat[Convergence behavior.]
    {
     \label{fig::Cube4HolesConvRefConv}
       \includegraphics[width=0.48\textwidth]{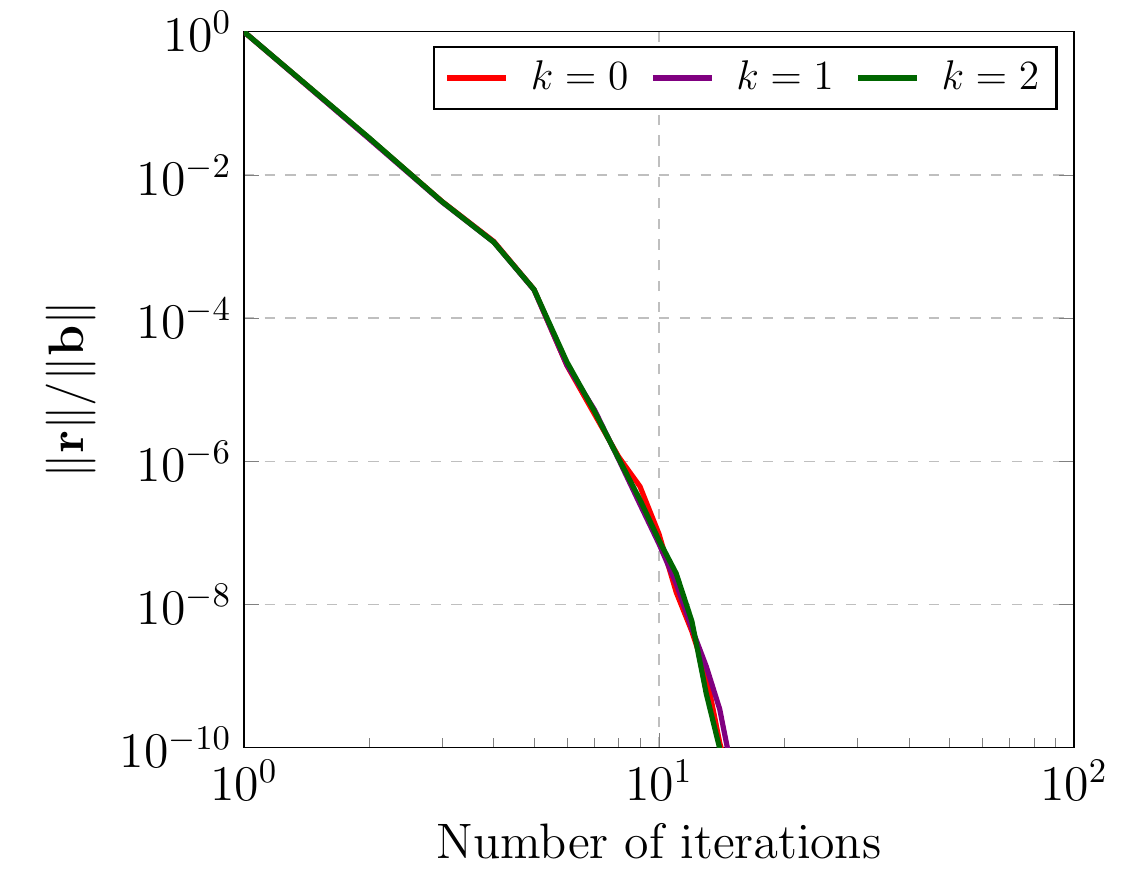}%
    }
    \subfloat[Computation cost of the hierarchical multigrid solver for a set of simulations running on 384 cores.]
    {
     \label{fig::Cube4HolesConvRefCost}
       \includegraphics[width=0.49\textwidth]{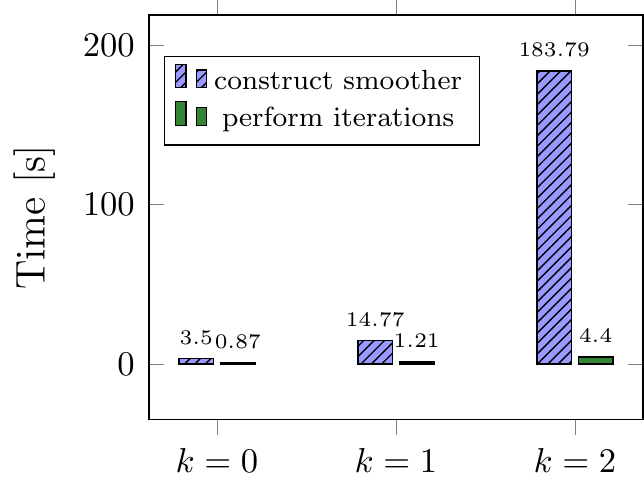}%
    }
 \caption{Performance of the CG solver with an $hp$-multigrid preconditioner and patchwise additive Schwarz smoothing for a three-dimensional benchmark involving multi-level $hp$-refinement.}
     \label{fig::Cube4HolesConvRef}
\end{figure}

The effect of the refinement level on the performance of the $hp$-multigrid preconditioner is now studied. Starting from an initial grid of $32^3$ elements with $p=2$, the mesh is refined in two steps towards the spherical cavities yielding a total of $4.1 \cdot 10^5$, $5.2 \cdot 10^5$ and $8.9 \cdot 10^5$ DOFs for $k=0$, $k=1$ and $k=2$, respectively. The patchwise additive Schwarz groups are used to construct the smoothers that are applied on every multigrid level $\ell \neq 0$ in analogy to the two-dimensional examples considered in the previous studies. This approach yields convergence rates that are independent of the refinement level as shown in Figure \ref{fig::Cube4HolesConvRefConv}.

The computational cost of the solver for the different refinement levels considered is shown in Figure \ref{fig::Cube4HolesConvRefCost}. The results shown are obtained in hybrid simulations on 394 cores that are partitioned into 32 MPI tasks each with 6 OMP threads. As expected, the inversion of the sub-matrices for the construction of the patchwise AS smoothers is the most time-consuming operation. This procedure is performed using the sparse direct solver \texttt{Pardiso}, since it outperforms the built-in invert function of the \texttt{BOOST} library \cite{Schling:2011:BCL:2049814} for large matrices. The time spent setting up the patchwise smoothers can be further reduced by the use of a more optimized inversion algorithm. Moreover, it is possible to conceive different strategies for the grouping of basis functions that yield smaller group sizes and therefore lower computation times.

\subsubsection{Loading of an aluminum rod}
\begin{figure}[H]
 \begin{center}
  \begin{minipage}[c]{0.52\textwidth}
   \begin{center}
    \includegraphics[width=0.9\textwidth]{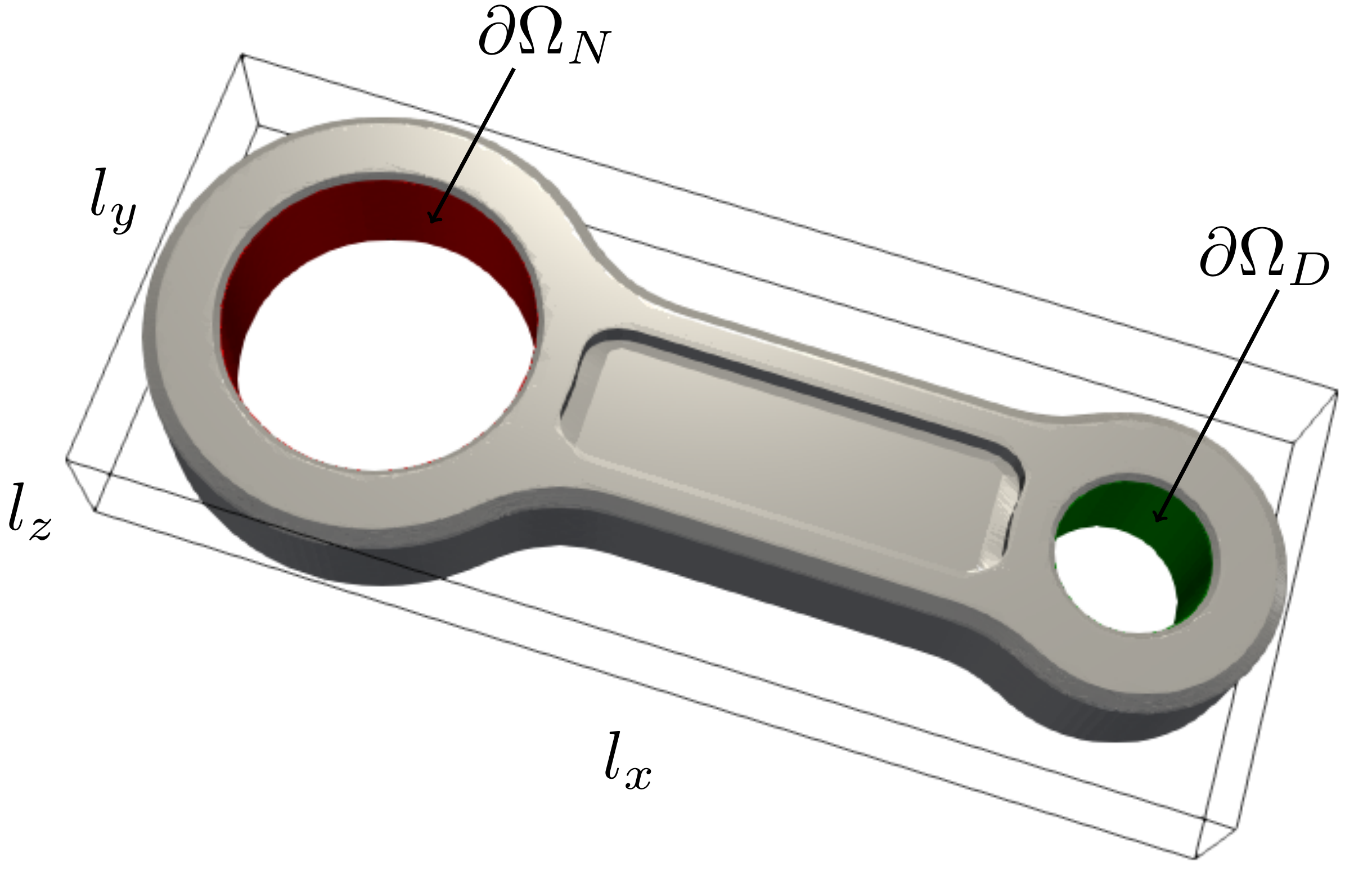}\newline%
       \small{(a) Rod geometry and boundary conditions}
   \end{center}
  \end{minipage}
 \begin{minipage}[c]{0.45\textwidth}
  \scriptsize
  \begin{center}
   \includegraphics[width=0.72\textwidth]{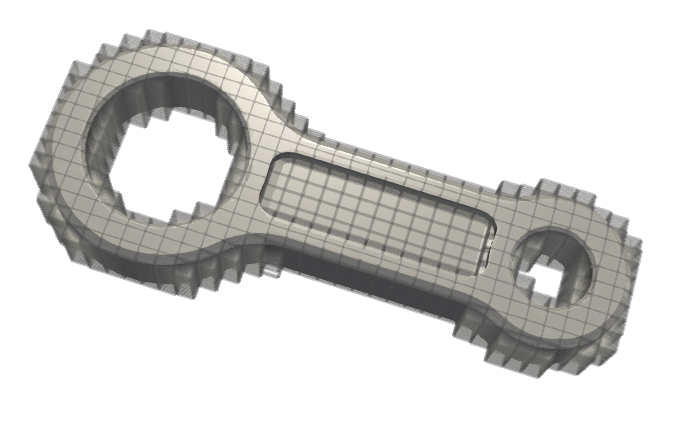}%
   \\ (b) Finite cell mesh.
   \\
   \vspace{3mm}
   \includegraphics[width=0.70\textwidth]{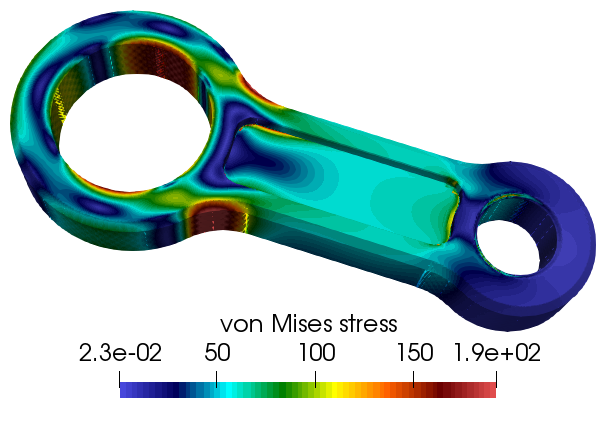}%
   \\ (c) Von Mises stress
  \end{center}
  \end{minipage}
  \caption{Loading of an aluminum connecting rod.}
  \label{fig::aluRod}
 \end{center}
\end{figure}

The next example considers the loading of an aluminum rod with an elastic modulus $\text{E} = 70\, \text{GPa}$ and a Poisson's ratio $\nu = 0.3$. Figure \hyperref[fig::aluRod]{\ref*{fig::aluRod}a} shows the rod's geometry with $l_x$ = 150 mm, $l_y$ = 60 mm and $l_z$ = 20 mm. The cylindrical surfaces labeled $\Omega_D$ are fixed using the penalty method with a penalty parameter $\beta = 10^6$. A surface load $t_x =10\, \text{N/mm}^2$ is applied on the surfaces labeled $\Omega_N$ in the direction of the rod's shaft. The resultant force acting on the rod is $F_x = 2\pi r h \cdot t_x = 2 \pi \cdot 20 \cdot 10 \approx 12.566\, \text{kN}$. Three finite cell discretizations with $h \in \lbrace{ 2, 1, 0.5 \rbrace}$ are considered and the polynomial degree of the grids is chosen as $p \in \lbrace{2 , 4 \rbrace}$. Table \ref{tab::aluminumrod} summarizes the number of DOFs in each grid. An octree scheme with a depth of 3 is used for the integration of the element matrices and the righthand side. 

\begin{table}[H]
    \centering
    \begin{tabularx}{0.8\textwidth}{|c|c| *{3}{Y|}}
\toprule
         &  & \multicolumn{3}{c|}{DOFs} \\
\toprule
        $h$ & mesh resolution & $p = 2$ & $p=3$ & $p=4$ \\
\toprule
        2 & $75 \times 30 \times 10$ & 156\,594 & 271\,596 & 492\,285 \\
        1 & $150 \times 60 \times 20$ & 1\,066\,059 & 1\,861\,425 & 3\,410\,673 \\
        0.5 & $300 \times 120 \times 40$ & 7\,880\,868 &  13\,755\,963 & 25\,365\,804 \\
\bottomrule
\end{tabularx}
        \caption{Summary of the number of DOFs for the different discretizations considered in the aluminum rod example.}
        \label{tab::aluminumrod}
\end{table}

\begin{table}[H]
\centering
\small
    \subfloat
{
\begin{tabularx}{0.5\textwidth}{|c| *{3}{Y|}}
\toprule
   & \multicolumn{3}{c|}{CG with elementwise AS}\\
	\cline{2-4}  & \multicolumn{3}{c|}{preconditioner}\\
\toprule
    \diagbox{$p$}{$h$} & 2 & 1 & 0.5 \\
\toprule
	2 & 766 &  1502 &  2966 \\
	3 & 770 &  1508 &  2973 \\
	4 & 835 &  1593 &  3254 \\
\bottomrule
\end{tabularx}
}
    \subfloat
{
\begin{tabularx}{0.5\textwidth}{|c| *{3}{Y|}}
\toprule
   & \multicolumn{3}{c|}{CG with multigrid preconditioner}\\
    \cline{2-4}  & \multicolumn{3}{c|}{elementwise AS smoother (5,5)}\\
\toprule
    \diagbox{$p$}{$h$} & 2 & 1 & 0.5 \\
\toprule
	2 & 16 &  15 &  15 \\
	3 & 14 &  13 &  13 \\
	4 & 18 &  17 &  17 \\
\bottomrule
\end{tabularx}
}
    \caption{Aluminium rod example: Convergence of a CG solver with two different preconditioners; \textit{i)} elementwise additive Schwarz preconditioner and \textit{ii)} $p$-multigrid preconditioner with elementwise additive Schwarz smoothing. The considered grids consist of trunk space hexahedral elements with different polynomial orders.}
        \label{fig::aluminumrodConvergenceP}
\end{table}

\paragraph{Comparison to different iterative solvers}\mbox{}\\
The study at hand compares the performance of the proposed multigrid solution scheme to alternative iterative approaches. In the first part of this study, the convergence behavior of the multigrid scheme with elementwise AS smoothing is compared to a CG solver that uses the additive Schwarz technique as a preconditioner. The results for the different finite cell discretizations considered in this example are presented in Table \ref{fig::aluminumrodConvergenceP}. Both preconditioners yield the expected convergence rates i.e.\ the elementwise additive Schwarz preconditioner yields rates proportional to $h^{-1}$, while the use of the $p$-multigrid preconditioner leads to convergence rates independent of the mesh size. Moreover, there is only a minimal difference in the convergence behavior of different polynomial orders for all preconditioners.

Next, we compare the convergence behavior and execution time of the presented multigrid approach to different solvers and preconditioners available in the \texttt{AztecOO} package of \texttt{Trilinos}. The solvers applied include: \textit{i)} a CG solver with diagonal scaling, denoted by cg-diag, \textit{ii)} a CG solver with elementwise additive Schwarz preconditioning, denoted by cg-eas, \textit{iii)} the solution approach presented in this contribution that consists of a CG solver with a $p$-multigrid preconditioner and elementwise additive Schwarz smoothing, denoted by cgmg-eas \textit{iv)} a CG solver with an algebraic multigrid solver based on smoothed aggregation, denoted by cg-amg and \textit{v)} a GMRES solver with an incomplete LU preconditioner, gmres-ilu. Note that the solvers cg-diag, cg-amg and grmes-ilu are included in the \texttt{AztecOO} package. These different solvers are used to solve the linear systems arising from the aluminium rod with $h=2$ and $p=3$. All simulations are carried out on 192 cores that are divided into 32 MPI tasks each with 6 OpenMP threads.  

\begin{figure}[H]
  \centering
  \subfloat[Relative residual.]
  {
    \includegraphics[width=0.6\textwidth]{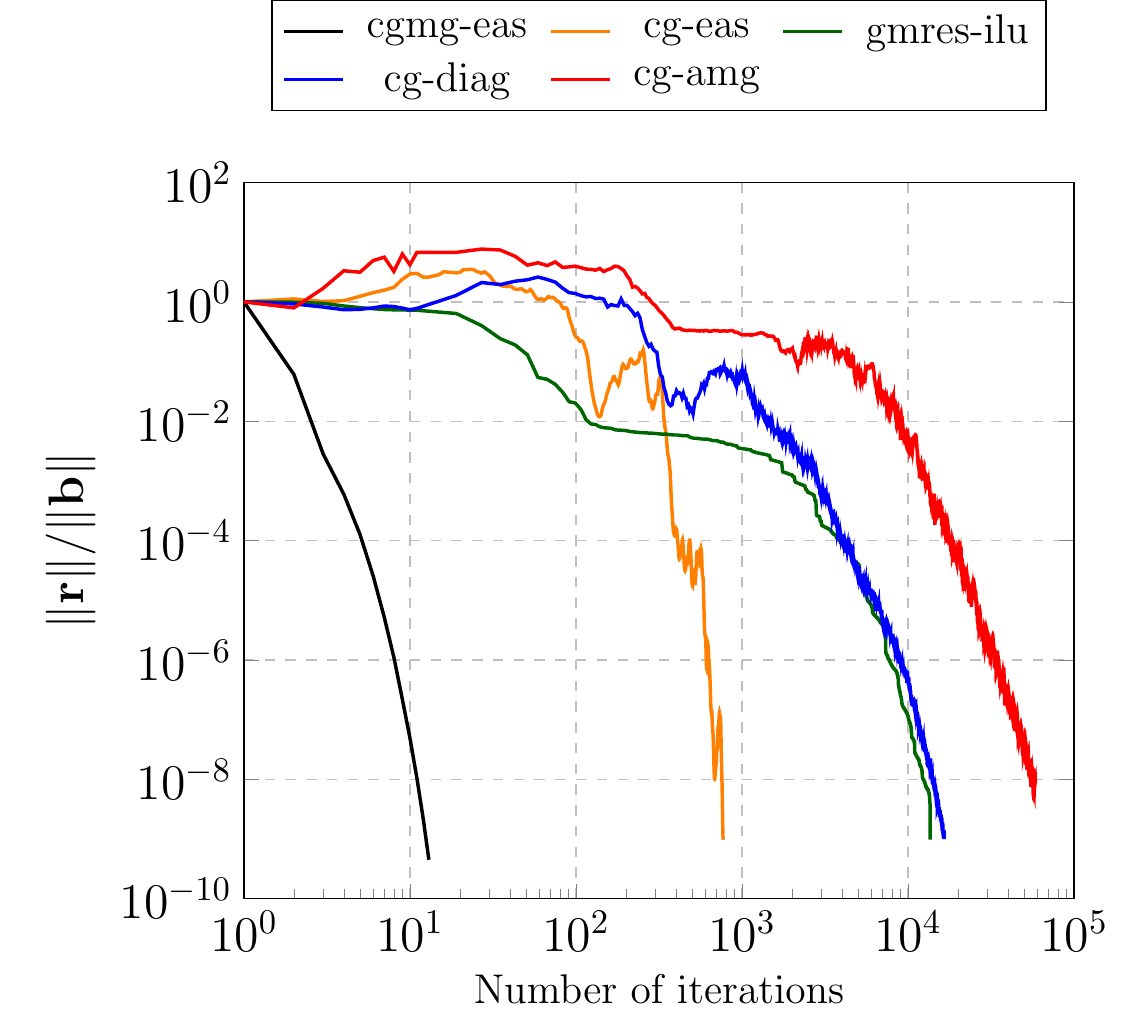}%
    \label{fig::solverComparisonConv}
  }
  \subfloat[Solver time.]
  {
    \includegraphics[width=0.4\textwidth]{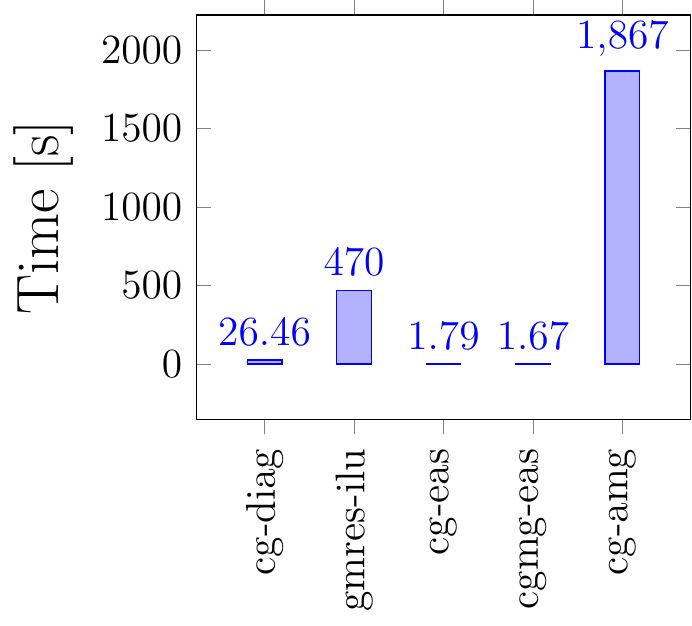}%
    \label{fig::solverComparisonTime}
  }
    \caption{Comparison of the convergence behavior and execution time of different iterative solvers in the aluminum rod example.}
        \label{fig::solverComparison}
\end{figure}

Figure \ref{fig::solverComparison} presents the comparison of the convergence behavior and execution time for the solution of a FCM problem on a complex geometry using five different solver and preconditioner configurations. The solver time shown in Figure \ref{fig::solverComparisonTime} includes both the construction and application of the solvers. The results affirm that the additive Schwarz-based solution techniques are well suited for large immersed systems as they not only have convergence rates superior to those of conventional solvers but also exhibit lower computational times. It should be noted that the poor performance of the cg-amg and gmres-ilu solvers is attributed to the ill-conditioning due to the cut elements, which both solvers fail to adequately address. 

\subsubsection{Weak scalability: Popcorn benchmark}
The final numerical example considers a classical immersed finite element benchmark on popcorn geometry. This geometry is commonly used in interface problems e.g.\ \citep{Chern2007,Chandrasekhar2012} and immersed analyses, see e.g. \citep{Burman2014}. A level-set function $\phi(x,y,z)$ is used to describe the surface of the popcorn geometry where

\begin{equation}
    \phi(x,y,z) = \sqrt{x^2 + y^2 + z^2} - r_0 - \sum \limits_{k=0}^{11} A e^{-((x - x _k)^2 +  (y - y _k)^2 + (z - z _k)^2 )/ \sigma^2},
\end{equation}
and
\begin{eqnarray}
    (x_k,y_k,z_k) &=& \frac{r_0}{\sqrt{5}} \Big( 2 \cos\Big(\frac{2k\pi}{5}\Big), \sin\Big(\frac{2k\pi}{5}\Big), 1 \Big), \hspace{3mm} \text{for} \hspace{3mm} k \in [0,4], \nonumber \\
    (x_k,y_k,z_k) &=& \frac{r_0}{\sqrt{5}} \Big( 2 \cos\Big(\frac{(2(k -5) - 1)\pi}{5}\Big), \sin\Big(\frac{(2(k -5) - 1)\pi}{5}\Big), 1 \Big), \hspace{3mm} \text{for} \hspace{3mm} k \in [5,9], \nonumber \\
    (x_{10},y_{10},z_{10}) &=& (0,0,r_0), \nonumber \\
    (x_{11},y_{11},z_{11}) &=& (0,0,-r_0), \nonumber
\end{eqnarray}
with the parameters $r_0$=0.6, $A$=3 and $\sigma$=0.2. The physical domain $\Omega_{\text{phys}}$ constitutes all points lying on the popcorn surface and its interior, i.e.\ all points with $\phi \leq 0$. Following \citep{Chern2007}, a Poisson problem with a prescribed solution field
 \begin{equation}
     u(x,y,z) = x^3 + xy^2 + y^3 + z^4 + \sin(3(x^2 + y^2)),  \hspace{3mm} \boldsymbol{x} \in
 \Omega_{phys}.
 \end{equation}
 is considered in the current example. This solution is defined through a source term $s=\Delta u $ and the enforcement of Dirichlet boundary conditions on $\phi = 0$. A marching cubes algorithm is run during the simulation to obtain the surfaces required for boundary condition application. The penalty method with $\beta = 10^4$ is used to apply the Dirichlet boundary conditions and the value of $\alpha$ is set to $10^{-6}$. An embedding domain consisting of a cube $(-1,1)^3$ is used in all numerical investigations.
 
\begin{figure}[H]
    \centering
    \subfloat[Popcorn geometry.]
    {
      \includegraphics[width=0.55\textwidth]{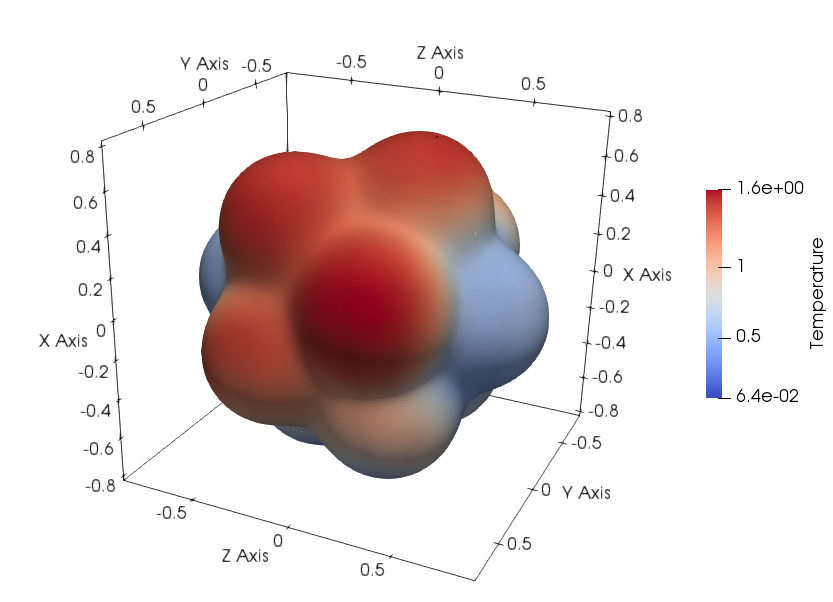}
      \label{fig::popcornGeometry}
    }
    \subfloat[Finite cell mesh.]
    {
      \includegraphics[width=0.55\textwidth]{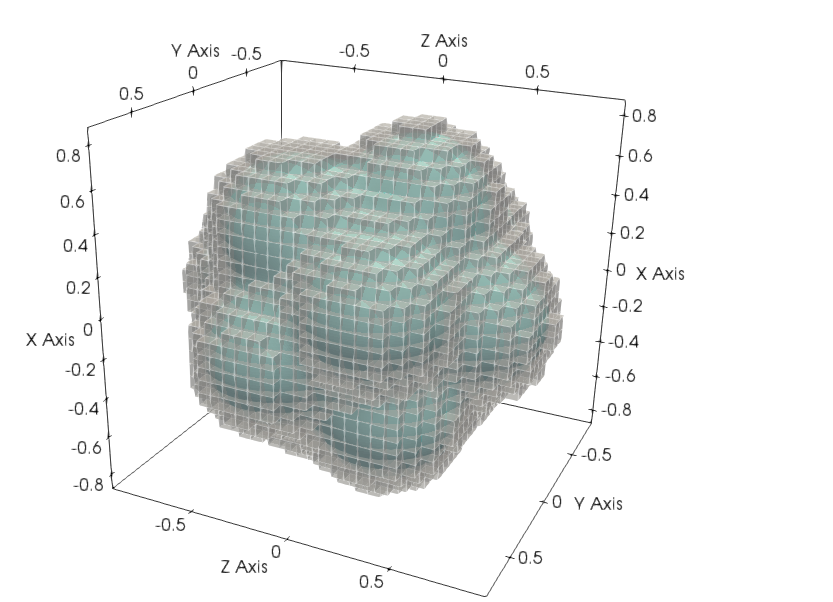}
    }
    \caption{Poisson problem posed on a popcorn domain.}
    \label{fig::popcornProblem}
\end{figure}
A weak scaling analysis is performed to investigate the efficiency of the proposed hierarchical multigrid approach. The simulations are performed on the SuperMUC-NG supercomputer at the Leibniz Supercomputing Center in Garching, Germany. The number of elements is increased in 12 steps while at the same time increasing the number of compute nodes $n_{\text{nodes}}$, from 1 to 2048. The 48 cores within a node are partitioned such that each node has a total of 8 MPI tasks and 6 OpenMP threads per task. The computational domain is generated in a fully parallel manner and discretized using hexahedral trunk space elements with a polynomial order $p=3$. The number of elements per MPI task is approximately 41\,000 in all simulations. The simulation with the coarsest mesh comprising a total of 250\,336 elements and $1.83 \cdot 10^6$ DOFs is run on 48 cores while the one with the finest mesh is run on 98\,304 cores, and has a total of 460\,980\,224 elements and approximately $3.2 \cdot 10^9$ DOFs. The different linear systems are solved using a parallel CG solver with a $p$-multigrid preconditioner and elementwise additive Schwarz smoothing.

\begin{figure}[H]
    \centering
    \subfloat[Convergence behavior.]
    {
      \includegraphics[width=0.42\textwidth]{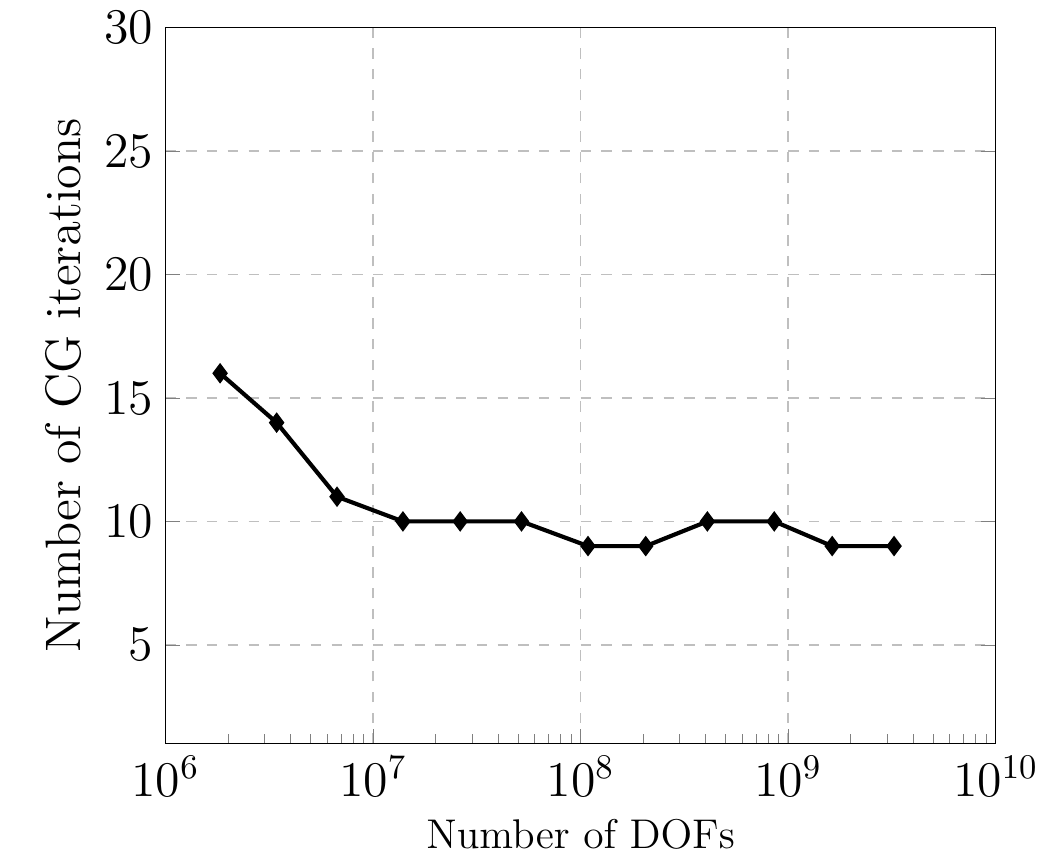}
    \label{fig::popcornIterations}
    }
    \subfloat[Solver time.]
    {
      \includegraphics[width=0.45\textwidth]{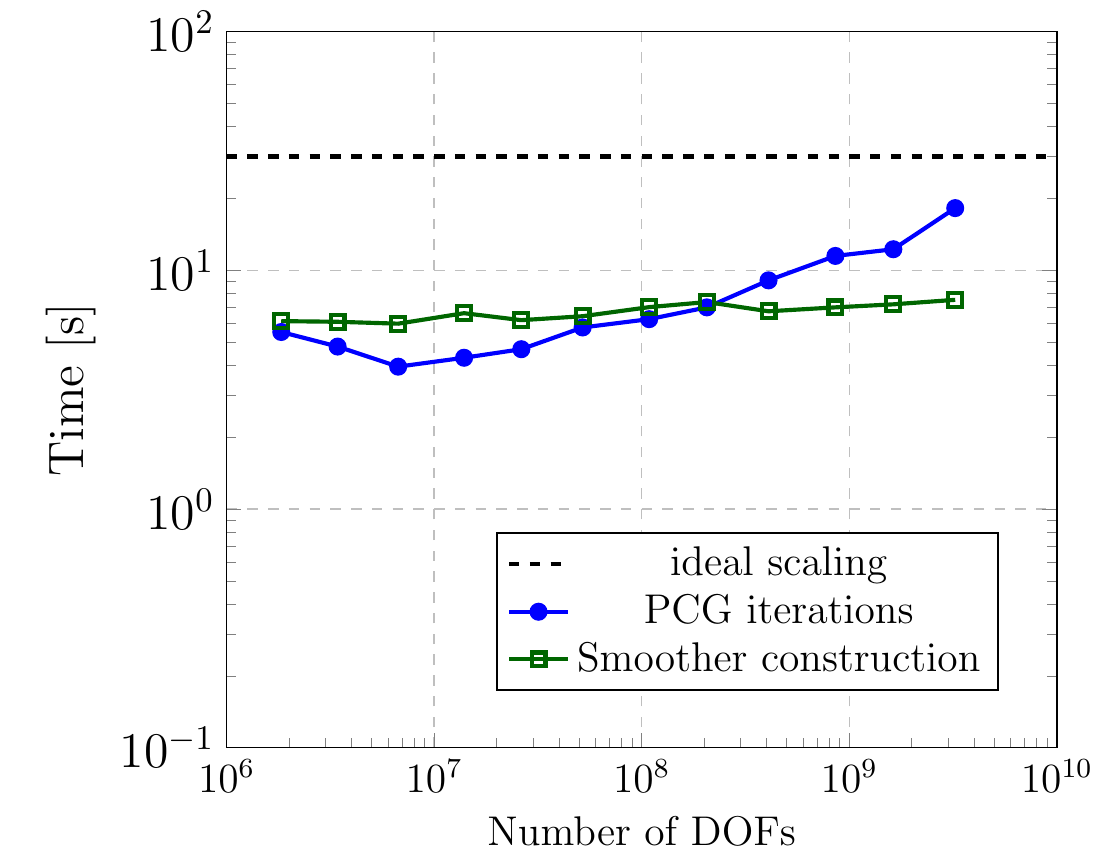}
    \label{fig::popcornTime}
    }
    \caption{Weak scaling analysis of the popcorn benchmark.}
    \label{fig::popcornExample}
\end{figure}

Figure \ref{fig::popcornExample} shows the results of the weak scaling analysis. It reports the number of iterations needed by the solver and the time spent to construct the additive Schwarz smoother and solve linear systems. The multigrid approach proposed in this manuscript shows favorable weak scaling, as the computational time does not significantly increase with an increase of the problem size.

\section{Conclusions}
The paper at hand presents a hierarchical multigrid method for immersed high-order discretizations based on the finite cell method and multi-level $hp$-refinement. This scheme robustly deals with ill-conditioning due to cut elements resulting in convergence rates independent of the cut configuration, element size and in certain scenarios, the polynomial order. The applicability of the scheme for immersed finite cell analyses on large distributed memory systems is portrayed in a collection of second-order problems arising from the Poisson equation and linear elasticity.

The cornerstone of this contribution is the development of a simple and efficient multigrid scheme that exploits the hierarchical nature of the basis functions in the finite cell method and the multi-level $hp$-scheme. For uniform grids, the multigrid levels are constructed using an arithmetic $p$-sequence that progressively reduces to polynomial order until $p = 1$. The nested subspace of a certain polynomial order is spanned by the basis functions up to and including this order. The $p$-multigrid approach is combined with an $h$-coarsening of the refined elements in the case of multi-level $hp$-refined grids. The main benefit of the proposed scheme is its simplicity, as the restriction and prolongation operations reduce to binary operations that can be easily performed without the assembly of restriction and prolongation operators.

Additive Schwarz-based smoothers are applied in this contribution to deal with the ill-conditioning due to cut cells. Two different procedures are suggested for finite cell discretizations: an elementwise approach is applied to uniform grids and a patchwise approach is used for multi-level $hp$-refined grids. Both smoothers are shown to sufficiently deal with small modes that occur on cut cells and yield convergence rates independent of the cut configuration and the mesh size $h$. The patchwise additive Schwarz smoother is shown to yield convergence rates independent of the polynomial order in the two dimensional numerical tests.

Furthermore, we demonstrate that the proposed multigrid scheme can be implemented within an efficient finite element framework in a massively parallel environment. A series of numerical examples demonstrate the suitability of the presented solution technique for large-scale immersed analysis. Faster convergence rates are achieved when the multigrid scheme is used as a preconditioner within a CG method rather than a stand-alone solver. The execution time for selected benchmark examples is also provided. These results indicated that the presented multigrid scheme exhibits low computation times as well as favorable weak scalability. The proposed scheme is shown to outperform standard iterative solvers both in terms of iteration counts and execution time for immersed systems.

Although the solution scheme presented in this contribution significantly improves the solution of high-order immersed systems involving $hp$-refinement, it can be further developed to increase its applicability. The construction of the patchwise preconditioner for high polynomials orders and several levels of refinement can be improved by applying more efficient matrix inversion techniques. This can be combined with the use of alternative strategies for the selection of additive Schwarz groups. These strategies can be designed to yield smaller group sizes and lower execution times. In this work, additive Schwarz smoothers are applied on every multigrid level with $\ell \neq 0$. The use of variable smoothing techniques on different multigrid levels is an important subject that can be addressed in future research. Apart from improving the efficiency of the solver, it would be also beneficial to invest time in a thorough mathematical analysis of the proposed multigrid method. This entails the derivation of bounds for the condition numbers as well as mathematical proofs for the multigrid method in an immersed setting.

\section*{Acknowledgements} 
The authors gratefully acknowledge the Competence Network for Scientific High Performance Computing in Bavaria (KONWIHR) and the Gauss Centre for Supercomputing e.V. (www.gauss-centre.eu) for the financial support and computing time provided on the SUPERMUG-NG at Leibniz Supercomputing Centre (www.lrz.de). We also extend our gratitude to the International Graduate School of Science and Engineering (IGSSE) of the Technical University of Munich for its financial support as well as the Deutsche Forschungsgemeinschaft (DFG, German Research Foundation) – Projektnummer 414265976 – TRR 277.

\clearpage
\bibliographystyle{ieeetr}
\bibliography{references.bib}

\begin{thebibliography}{10}

\bibitem{Parvizian2007}
J.~Parvizian, A.~D{\"u}ster, and E.~Rank, ``Finite cell method,'' {\em
  Computational Mechanics}, vol.~41, no.~1, pp.~121--133, 2007.

\bibitem{Duster2008}
A.~D{\"u}ster, J.~Parvizian, Z.~Yang, and E.~Rank, ``The finite cell method for
  three-dimensional problems of solid mechanics,'' {\em Computer Methods in
  Applied Mechanics and Engineering}, vol.~197, no.~45\textendash{}48,
  pp.~3768--3782, 2008.

\bibitem{Burman2012}
E.~Burman and P.~Hansbo, ``Fictitious domain finite element methods using cut
  elements: {{II. A stabilized Nitsche method}},'' {\em Applied Numerical
  Mathematics}, vol.~62, no.~4, pp.~328 -- 341, 2012.

\bibitem{Burman2014}
E.~Burman, S.~Claus, P.~Hansbo, M.~G. Larson, and A.~Massing, ``{{CutFEM:}}
  {{Discretizing}} geometry and partial differential equations,'' {\em
  International Journal for Numerical Methods in Engineering}, vol.~104, no.~7,
  pp.~472--501, 2014.

\bibitem{Badia2018a}
S.~Badia, F.~Verdugo, and A.~F. Mart{\'{\i}}n, ``The aggregated unfitted finite
  element method for elliptic problems,'' {\em Computer Methods in Applied
  Mechanics and Engineering}, vol.~336, pp.~533--553, 2018.

\bibitem{Soriano2013}
E.~Nadal~Soriano, J.~R{\'o}denas, J.~Albelda, M.~Tur, J.~Taranc{\'o}n, and
  F.~Fuenmayor, ``Efficient finite element methodology based on cartesian
  grids: Application to structural shape optimization,'' {\em Abstract and
  Applied Analysis}, 2013.

\bibitem{Navarro2018}
J.~M. Navarro-Jim{\'e}nez, M.~Tur, F.~J. Fuenmayor, and J.~J. R{\'o}denas, ``On
  the effect of the contact surface definition in the cartesian grid finite
  element method,'' {\em Advanced Modeling and Simulation in Engineering
  Sciences}, vol.~5, p.~12, May 2018.

\bibitem{main2018}
A.~Main and G.~Scovazzi, ``{The shifted boundary method for embedded domain
  computations. Part I: Poisson and Stokes problems},'' {\em Journal of
  Computational Physics}, vol.~372, pp.~972--995, 2018.

\bibitem{main2018b}
A.~Main and G.~Scovazzi, ``{The shifted boundary method for embedded domain
  computations. Part II: Linear advection-diffusion and incompressible
  Navier-Stokes equations},'' {\em Journal of Computational Physics}, vol.~372,
  pp.~996--1026, 2018.

\bibitem{Prenter2017}
F.~de~Prenter, {\relax C.V}.~Verhoosel, {\relax G.J}.~van Zwieten, and {\relax
  E.H}.~van Brummelen, ``Condition number analysis and preconditioning of the
  finite cell method,'' {\em Computer Methods in Applied Mechanics and
  Engineering}, vol.~316, pp.~297--327, 2017.

\bibitem{Badia2017}
S.~Badia and F.~Verdugo, ``Robust and scalable domain decomposition solvers for
  unfitted finite element methods,'' {\em Journal of Computational and Applied
  Mathematics}, vol.~344, pp.~740--759, 2017.

\bibitem{Prenter2019}
F.~de~Prenter, C.~Verhoosel, and E.~van Brummelen, ``Preconditioning immersed
  isogeometric finite element methods with application to flow problems,'' {\em
  Computer Methods in Applied Mechanics and Engineering}, vol.~348,
  pp.~604--631, 2019.

\bibitem{Jomo2019}
J.~Jomo, F.~de~Prenter, M.~Elhaddad, D.~D'Angella, C.~Verhoosel,
  S.~Kollmannsberger, J.~Kirschke, V.~N{\"u}bel, E.~van Brummelen, and E.~Rank,
  ``Robust and parallel scalable iterative solutions for large-scale finite
  cell analyses,'' {\em Finite Elements in Analysis and Design}, vol.~163,
  pp.~14--30, 2019.

\bibitem{Burman2010}
E.~Burman, ``Ghost penalty,'' {\em Comptes Rendus Mathematique}, vol.~348,
  no.~21, pp.~1217 -- 1220, 2010.

\bibitem{Elfverson2018}
D.~Elfverson, M.~G. Larson, and K.~Larsson, ``{{CutIGA with basis function
  removal}},'' {\em Advanced Modeling and Simulation in Engineering Sciences},
  vol.~5, no.~1, p.~6, 2018.

\bibitem{Marussig2017review}
B.~Marussig and T.~Hughes, ``A {R}eview of {T}rimming in {I}sogeometric
  {A}nalysis: {C}hallenges, {D}ata {E}xchange and {S}imulation {A}spects,''
  {\em Archives of Computational Methods in Engineering}, vol.~25, no.~4,
  pp.~1059--1127, 2017.

\bibitem{Hackbusch1982}
W.~Hackbusch and U.~Trottenberg, {\em {Multigrid Methods. Proceedings of the
  Conference Held at K{\"o}ln-Porz, November 23-27, 1981}}.
\newblock 01 1982.

\bibitem{hackbusch2013}
W.~Hackbusch, {\em {Multi-Grid Methods and Applications}}.
\newblock Springer Series in Computational Mathematics, Springer Berlin
  Heidelberg, 2013.

\bibitem{trottenberg2001}
U.~Trottenberg, C.~Ulrich~Trottenberg, C.~Oosterlee, A.~Schuller, A.~Brandt,
  P.~Oswald, and K.~St{\"u}ben, {\em Multigrid}.
\newblock Elsevier Science, 2001.

\bibitem{Briggs2000}
W.~Briggs, V.~Henson, and S.~McCormick, {\em A Multigrid Tutorial, 2nd
  Edition}.
\newblock 01 2000.

\bibitem{wesseling2004}
P.~Wesseling, {\em {An Introduction to Multigrid Methods}}.
\newblock An Introduction to Multigrid Methods, R.T. Edwards, 2004.

\bibitem{shapira2003}
Y.~Shapira, {\em {Matrix-Based Multigrid: Theory and Applications}}.
\newblock Numerical methods and algorithms, Kluwer Academic Publishers, 2003.

\bibitem{Graig1985}
A.~W. Craig and O.~C. Zienkiewicz, ``A multigrid algorithm using a hierarchical
  finite element basis,'' in {\em Multigrid Methods for Integral and
  Differential Equations} (D.~J. Paddon and Holstein, eds.), pp.~310--312,
  Oxford,: Clarendon Press, 1985.

\bibitem{Babuska1989}
I.~Babu\v{s}ka, M.~Griebel, and J.~Pitk{\"a}ranta, ``The problem of selecting
  the shape functions for a p-type finite element,'' {\em International Journal
  for Numerical Methods in Engineering}, vol.~28, no.~8, pp.~1891--1908, 1989.

\bibitem{Hofreither2016}
C.~Hofreither, B.~J{\"u}ttler, G.~Kiss, and W.~Zulehner, ``{Multigrid methods
  for isogeometric analysis with THB-splines},'' {\em Computer Methods in
  Applied Mechanics and Engineering}, vol.~308, pp.~96--112, 2016.

\bibitem{Gahalaut2013}
K.~Gahalaut, J.~Kraus, and S.~Tomar, ``Multigrid methods for isogeometric
  discretization,'' {\em Computer Methods in Applied Mechanics and
  Engineering}, vol.~253, pp.~413--425, 2013.

\bibitem{Tielen2020}
R.~Tielen, M.~M\"{o}ller, D.~G\"{o}ddeke, and C.~Vuik, ``{$p$-multigrid methods
  and their comparison to $h$-multigrid methods within Isogeometric
  Analysis},'' {\em Computer Methods in Applied Mechanics and Engineering},
  vol.~372, p.~113347, 2020.

\bibitem{Riva2019}
A.~P. de~la Riva, C.~Rodrigo, and F.~J. Gaspar, ``{A Robust Multigrid Solver
  for Isogeometric Analysis Based on Multiplicative Schwarz Smoothers},'' {\em
  SIAM Journal on Scientific Computing}, vol.~41, pp.~321--345, 2019.

\bibitem{Donatelli2015}
M.~Donatelli, C.~Garoni, C.~Manni, S.~Serra-Capizzano, and H.~Speleers,
  ``Robust and optimal multi-iterative techniques for {I}g{A} {G}alerkin linear
  systems,'' {\em Computer Methods in Applied Mechanics and Engineering},
  vol.~284, pp.~230--264, 2015.

\bibitem{Donatelli2017}
M.~Donatelli, C.~Garoni, C.~Manni, S.~Serra-Capizzano, and H.~Speleers,
  ``{S}ymbol-{B}ased {M}ultigrid {M}ethods for {G}alerkin {B}-{S}pline
  {I}sogeometric {A}nalysis,'' {\em SIAM Journal on Numerical Analysis},
  vol.~55, no.~1, pp.~31--62, 2017.

\bibitem{Bracco2019}
C.~Bracco, D.~Cho, C.~Giannelli, and R.~Vazquez, ``{BPX} {P}reconditioners for
  {I}sogeometric {A}nalysis {U}sing ({T}runcated) {H}ierarchical {B}-splines,''
  2019.

\bibitem{Verdugo2019}
F.~Verdugo, A.~F. Mart{\'{\i}}n, and S.~Badia, ``Distributed-memory
  parallelization of the aggregated unfitted finite element method,'' {\em
  Computer Methods in Applied Mechanics and Engineering}, vol.~357, p.~112583,
  2019.

\bibitem{badia2020aggregated}
S.~Badia, A.~F. Mart{\'{\i}}n, E.~Neiva, and F.~Verdugo, ``The aggregated
  unfitted finite element method on parallel tree-based adaptive meshes,'' {\em
  preprint available on arXiv}, 2020.

\bibitem{badia2020generic}
S.~Badia, A.~F. Mart{\'{\i}}n, E.~Neiva, and F.~Verdugo, ``A generic finite
  element framework on parallel tree-based adaptive meshes,'' {\em preprint
  available on arXiv}, 2020.

\bibitem{Nussig2018}
A.~N\"u{\ss}ing, {\em {Fitted and unitted finite element methods for solving
  the EEG forward problem}}.
\newblock {PhD Thesis}, University of M{\"u}nster, M\"unster, 2018.

\bibitem{Prenter2019b}
F.~de~Prenter, C.~V. Verhoosel, E.~H. van Brummelen, J.~A. Evans, C.~Messe,
  J.~Benzaken, and K.~Maute, ``Multigrid solvers for immersed finite element
  methods and immersed isogeometric analysis,'' {\em Computational Mechanics},
  2019.

\bibitem{Saberi2020}
S.~Saberi, A.~Vogel, and G.~Meschke, ``Parallel finite cell method with
  adaptive geometric multigrid,'' in {\em Euro-Par 2020: Parallel Processing}
  (M.~Malawski and K.~Rzadca, eds.), (Cham), pp.~578--593, Springer
  International Publishing, 2020.

\bibitem{Roquist1987}
E.~M. R{\o}nquist and A.~Patera, ``{Spectral element multigrid. I. Formulation
  and numerical results},'' {\em Journal of Scientific Computing}, vol.~2,
  pp.~389--406, 1987.

\bibitem{Yserentant1985}
H.~Yserentant, ``Hierarchical bases of finite-element spaces in the
  discretization of nonsymmetric elliptic boundary value problems,'' {\em
  Computing}, vol.~35, pp.~39--49, 1985.

\bibitem{Yserentant1986}
H.~Yserentant, ``Hierarchical bases give conjugate gradient type methods a
  multigrid speed of convergence,'' {\em Applied Mathematics and Computation},
  vol.~19, no.~1, pp.~347--358, 1986.

\bibitem{Foresti1989}
S.~Foresti, G.~Brussino, S.~Hassanzadeh, and V.~Sonnad, ``Multilevel solution
  of the p-version of finite elements,'' {\em Computer Physics Communications},
  vol.~53, pp.~349--355, 05 1989.

\bibitem{Mitchell2010}
W.~F. Mitchell, ``The hp-multigrid method applied to hp-adaptive refinement of
  triangular grids,'' {\em Numerical Linear Algebra with Applications},
  vol.~17, no.~2-3, pp.~211--228, 2010.

\bibitem{Nastase2006}
C.~R. Nastase and D.~J. Mavriplis, ``{High-order discontinuous Galerkin methods
  using an hp-multigrid approach},'' {\em Journal of Computational Physics},
  vol.~213, no.~1, pp.~330 -- 357, 2006.

\bibitem{Yserentant1986b}
H.~Yserentant, ``On the multi-level splitting of finite element spaces,'' {\em
  Numerische Mathematik}, vol.~49, pp.~379--412, 1986.

\bibitem{Zander2015}
N.~Zander, T.~Bog, S.~Kollmannsberger, D.~Schillinger, and E.~Rank,
  ``{{Multi-Level hp-Adaptivity: High-Order Mesh Adaptivity without the
  Difficulties of Constraining Hanging Nodes}},'' {\em Computational
  Mechanics}, vol.~55, no.~3, pp.~499--517, 2015.

\bibitem{Zander2016}
N.~Zander, T.~Bog, M.~Elhaddad, F.~Frischmann, S.~Kollmannsberger, and E.~Rank,
  ``{{The Multi-Level hp-Method for Three-Dimensional Problems: {{Dynamically}}
  Changing High-Order Mesh Refinement with Arbitrary Hanging Nodes}},'' {\em
  Computer Methods in Applied Mechanics and Engineering}, vol.~310,
  pp.~252--277, 2016.

\bibitem{Babuska1981}
I.~Babuska, B.~Szabo, and I.~Katz, ``The p-{Version} of the {Finite} {Element}
  {Method},'' {\em SIAM Journal on Numerical Analysis}, vol.~18, no.~3,
  pp.~515--545, 1981.

\bibitem{Elhaddad2017}
M.~Elhaddad, N.~Zander, T.~Bog, L.~Kudela, S.~Kollmannsberger, J.~Kirschke,
  T.~Baum, M.~Ruess, and E.~Rank, ``Multi-level hp-finite cell method for
  embedded interface problems with application in biomechanics,'' {\em
  International Journal for Numerical Methods in Biomedical Engineering},
  vol.~34, no.~4, p.~e2951, 2018.

\bibitem{Ozcan2018}
A.~{\"O}zcan, S.~Kollmannsberger, J.~Jomo, and E.~Rank, ``Residual stresses in
  metal deposition modeling: Discretizations of higher order,'' {\em Computers
  \& Mathematics with Applications}, 2018.

\bibitem{Hug2020}
L.~Hug, S.~Kollmannsberger, Z.~Yosibash, and E.~Rank, ``{A 3D benchmark problem
  for crack propagation in brittle fracture},'' {\em Computer Methods in
  Applied Mechanics and Engineering}, vol.~364, p.~112905, 2020.

\bibitem{Babuska1973}
I.~Babu\v{s}ka, ``The {{Finite Element Method}} with {{Penalty}},'' {\em
  Mathematics of Computation}, vol.~27, no.~122, p.~221, 1973.

\bibitem{Schillinger2014}
D.~Schillinger and M.~Ruess, ``The {{Finite Cell Method}}: {{A Review}} in the
  {{Context}} of {{Higher}}-{{Order Structural Analysis}} of {{CAD}} and
  {{Image}}-{{Based Geometric Models}},'' {\em Archives of Computational
  Methods in Engineering}, vol.~22, no.~3, pp.~391--455, 2014.

\bibitem{Duster2017}
A.~D{\"u}ster, E.~Rank, and B.~A. Szab{\'o}, ``The p-version of the finite
  element method and finite cell methods,'' in {\em Encyclopedia of
  {{Computational}} Mechanics}, vol.~2, pp.~1--35, Chichester, West Sussex:
  {John Wiley \& Sons}, 2017.

\bibitem{DAngella2016}
D.~D'Angella, N.~Zander, S.~Kollmannsberger, F.~Frischmann, E.~Rank,
  A.~Schr{\"o}der, and A.~Reali, ``Multi-level hp-adaptivity and explicit error
  estimation,'' {\em Advanced Modeling and Simulation in Engineering Sciences},
  vol.~3, no.~1, p.~33, 2016.

\bibitem{Darrigrand2020}
V.~Darrigrand, D.~Pardo, T.~Chaumont-Frelet, I.~G{\'o}mez-Revuelto, and L.~E.
  Garcia-Castillo, ``A painless automatic hp-adaptive strategy for elliptic
  problems,'' {\em Finite Elements in Analysis and Design}, vol.~178,
  p.~103424, 2020.

\bibitem{Zander2016b}
N.~Zander, {\em Multi-Level hp-{{FEM}}: Dynamically Changing High-Order Mesh
  Refinement with Arbitrary Hanging Nodes}.
\newblock {PhD Thesis}, Technische Universit{\"a}t M{\"u}nchen, Munich, 2016.

\bibitem{Saad2003}
Y.~Saad, {\em Iterative Methods for Sparse Linear Systems}.
\newblock Philadelphia, PA, USA: Society for Industrial and Applied
  Mathematics, 2nd~ed., 2003.

\bibitem{mpi1994}
M.~P. Forum, ``{MPI: A Message-Passing Interface Standard},'' tech. rep.,
  Knoxville, TN, USA, 1994.

\bibitem{openmp08}
{OpenMP Architecture Review Board}, ``{{OpenMP} Application Program Interface
  Version 3.0},'' 2008.

\bibitem{Trilinos}
M.~A. Heroux, R.~A. Bartlett, V.~E. Howle, R.~J. Hoekstra, J.~J. Hu, T.~G.
  Kolda, R.~B. Lehoucq, K.~R. Long, R.~P. Pawlowski, E.~T. Phipps, A.~G.
  Salinger, H.~K. Thornquist, R.~S. Tuminaro, J.~M. Willenbring, A.~Williams,
  and K.~S. Stanley, ``{An overview of the Trilinos project},'' {\em ACM
  Transactions on Mathematical Software}, vol.~31, no.~3, pp.~397--423, 2005.

\bibitem{Schenk2011}
O.~Schenk and K.~G{\"a}rtner, {\em PARDISO}, pp.~1458--1464.
\newblock Boston, MA: Springer US, 2011.

\bibitem{Schling:2011:BCL:2049814}
B.~Schling, {\em The Boost C++ Libraries}.
\newblock XML Press, 2011.

\bibitem{Chern2007}
I.-L. Chern and Y.-C. Shu, ``A coupling interface method for elliptic interface
  problems,'' {\em Journal of Computational Physics}, vol.~225, no.~2,
  pp.~2138--2174, 2007.

\bibitem{Chandrasekhar2012}
C.~Annavarapu, M.~Hautefeuille, and J.~Dolbow, ``{A robust Nitsche's
  formulation for interface problems},'' {\em Computer Methods in Applied
  Mechanics and Engineering}, vol.~s 225--228, pp.~44--54, 06 2012.

\end{thebibliography}

\end{document}